\documentclass[12pt]{amsart}
\usepackage{bbding}
\usepackage{dsfont}
\usepackage{graphicx}
\usepackage{amsthm,amscd}
\usepackage{calc}
\usepackage{hyperref}
\usepackage{mathrsfs}
\usepackage{mathrsfs}
\usepackage{CJK,fancyhdr,amscd}
\usepackage{anysize}
\usepackage{amsmath,amssymb,amsfonts}
\usepackage{epsfig}
\usepackage{latexsym}
\usepackage{hyperref}
\usepackage{indentfirst, latexsym}
\usepackage{graphics}
\usepackage{amsmath, bm}
\usepackage{amsfonts}
\usepackage{amssymb}%
\usepackage{color}
\usepackage{mathtools}
\providecommand{\U}[1]{\protect\rule{.1in}{.1in}}

\numberwithin{equation}{section}
\providecommand{\U}[1]{\protect\rule{.1in}{.1in}}

\newtheorem{theorem} {Theorem} [section]
\newtheorem{proposition}[theorem]{Proposition}
\newtheorem{corollary}  [theorem]     {Corollary}
\newtheorem{lemma}  [theorem]     {Lemma}
\newtheorem{example}  [theorem]     {Example}
\newtheorem{remark}  [theorem]     {Remark}
\newtheorem{definition}  [theorem]     {Definition}

\newcommand{\dmath}{\mathrm{d}\:\!}
\newcommand{\p}{\partial}
\newcommand{\bp}{\bar\partial}

\newcommand{\bthm}{\begin{theorem}}
\newcommand{\ethm}{\end{theorem}}
\newcommand{\blem}{\begin{lemma}}
\newcommand{\elem}{\end{lemma}}
\newcommand{\bprop}{\begin{proposition}}
\newcommand{\eprop}{\end{proposition}}
\newcommand{\brmk}{\begin{remark}}
\newcommand{\ermk}{\end{remark}}
\newcommand{\bcor}{\begin{corollary}}
\newcommand{\ecor}{\end{corollary}}
\newcommand{\bdefi}{\begin{definition}}
\newcommand{\edefi}{\end{definition}}
\newcommand{\beq}{\begin{equation}}
\newcommand{\eeq}{\end{equation}}
\newcommand{\bex}{\begin{example}}
\newcommand{\eex}{\end{example}}

\newcommand{\tOmega}{\widetilde{\Omega}}
\newcommand{\bOmega}{\overline{\Omega}{\mathstrut}}

\usepackage[lmargin=1in, rmargin=1in]{geometry}
\title{The Monge-Amp\`{e}re equation for $(n-1)$\nobreakdash-quaternionic PSH functions on a hyperK\"{a}hler manifold}
\author{Jixiang Fu \and Xin Xu \and Dekai Zhang}
\address{Shanghai Center for Mathematical Sciences, Jiangwan Campus, Fudan University, Shanghai, 200438, China}
\email{majxfu@fudan.edu.cn}
\address{School of Mathematical Sciences, Fudan University, Shanghai 200433, China}
\email{20110180015@fudan.edu.cn}
\address{Department of Mathematics, Shanghai University, Shanghai, 200444, China}
\email{dkzhang@shu.edu.cn}

\begin{document}
\maketitle
\begin{abstract}
    We prove the existence of unique smooth solutions to the quaternionic Monge-Amp\`{e}re equation for $(n-1)$\nobreakdash-quaternionic plurisubharmonic functions on a hyperK\"{a}hler manifold and thus obtain solutions for the quaternionic form type equation. We derive the $C^0$ estimate by establishing a Cherrier-type inequality as in Tosatti and Weinkove \cite{tosatti2017monge}. By adopting the approach of Dinew and Sroka \cite{dinew2021hkt} to our context, we obtain $C^1$ and $C^2$ estimates without assuming the flatness of underlying hyperK\"{a}hler metric comparing to previous results \cite{gentili2022fully}.
\end{abstract}
\section{Introduction}
\label{introduction}
A hypercomplex manifold is a smooth manifold $M$ together with a triple $(I, J, K)$ of
complex structures satisfying the quaternoinic relation
\[IJ = -JI = K.\]

Let $(M, I, J, K)$ be a hypercomplex manifold, and $g$ a Riemannian metric on $M$. 
The metric $g$ is called hyperhermitian if $g$ is hermitian with respect to $I, J, K$, i.e. for any vector fields $X$ and $Y$ on $M$,
\[g(X, Y) = g(XI, YI) = g(XJ, YJ) = g(XK, YK).\]
Denote by $\omega_I, \omega_J, \omega_K$ the fundamental form corresponding to $I, J, K$ respectively and  let $\Omega = \omega_J + i \omega_K$. On a hyperhermitian manifold $(M, I, J, K, g)$, the metric $g$ is called hyperK\"{a}hler (HK) if $d\Omega = 0$ 
or equivalently $d\omega_I = d\omega_J = d\omega_K = 0$, and called hyperK\"{a}hler with torsion (HKT) if $\partial\Omega = 0$.

In analogy with the classical Calabi-Yau theorem \cite{yau1978ricci} on the complex Monge-Amp\`{e}re equation on a K\"{a}hler manifold, Alesker and Verbitsky \cite{alesker2010quaternionic} conjectured the existsnce of solutions to the quaternoinic Monge-Amp\`{e}re equation on a compact HKT manifold of quaternionic dimension $n$. It takes the form
\begin{equation}
\label{eqn:intro-q-calabi}
    \begin{split}
    (\Omega + \partial \partial_J u)^n &= e^f \Omega^n, \\
    \Omega + \partial \partial_J u &> 0,
    \end{split}
\end{equation}
where $\partial_J = J^{-1} \circ \overline{\partial} \circ J $. While general solution to this equation remains open, partial results can be found in \cite{alesker2013solvability, BGV22, gentili2022quaternionic, dinew2021hkt, alesker2010quaternionic, ALESKER20171, alesker2006plurisubharmonic, sroka2020}. Specifically, Alesker and Verbitsky \cite{alesker2010quaternionic} obtained $C^0$ estimate when the canonical bundle is holomorphically trivial. Alesker \cite{alesker2013solvability} proved the conjecture on compact manifolds with a flat hyperK\"{a}hler metric. In \cite{ALESKER20171} Alesker and Shelukhin proved $C^0$ estimate without any extra assumptions and the proof was later simplyfied by Sroka \cite{sroka2020}. Dinew and Sroka \cite{dinew2021hkt} solved equation \eqref{eqn:intro-q-calabi} on a hyperK\"{a}hler manifold.

As in the complex setting, we consider the quaternoinic form-type equation, as the analogue of the form-type equation which was proposed and also solved on a K\"{a}hler manifold of nonnegative bisectional curvature by Fu-Wang-Wu \cite{fu2010form, fu2015form}. It was later shown by Tosatti and Weinkove \cite{tosatti2017monge} that the assumption on curvature can be removed. 

In particular, one can define quaternionic balanced metrics on hypercomplex manifolds 
by $\p \Omega^{n-1} = 0$ (see \cite{MR3645310}). Let $(M, I, J, K, g, \Omega)$ be a hypercomplex manifold of quaternionic dimension $n$, and $g_0$ a quaternionic balanced metric on $M$ with induced $(2, 0)$\nobreakdash-form $\Omega_0$.
Let $\varphi$ be a $(2n-4,0)$\nobreakdash-form such that $\Omega_0^{n-1} + \p \p_J \varphi$ is strictly positive. Then there exists a quaternionic balanced metric $\Omega_{\varphi}$ such that 
\begin{equation}
    \Omega_{\varphi}^{n-1} = \Omega_0^{n-1} + \p \p_J \varphi.
\end{equation}
The quaternionic form-type Calabi-Yau equation is written as
\begin{equation}
\label{eqn:form-type}
    \Omega_{\varphi}^n = e^{f' + b'}\Omega^n
\end{equation}
where $f'$ is a given smooth function on $M$ and $b'$ is a uniquely determined constant. Solving equation \eqref{eqn:form-type} gives a quaternionic balanced metric $\Omega_{\varphi}$ with prescribed volume form up to scaling. 
One can reduce the form-type equation to function type by considering a function $u \in C^{\infty}(M, \mathbb{R})$ such that
$\Omega_0^{n-1} + \partial \partial_J (u \Omega^{n-2})$ is strictly positive, and denote by
$\Omega_u$ the unique strictly positive $(2, 0)$\nobreakdash-form such that
\begin{equation}
\label{eqn:(n-1)power-1}
    \Omega_u^{n-1} = \Omega_0^{n-1} + \partial \partial_J (u \Omega^{n-2}).
\end{equation}
Then equation \eqref{eqn:form-type} is reduced to
\begin{equation}
\label{eqn:func-type}
    \Omega_u^n = e^{f' + b'}\Omega^n.
\end{equation}
In particular when $\Omega$ is HKT, $\p \Omega = \p_J \Omega = 0$. Then \eqref{eqn:(n-1)power-1} becomes
\begin{equation}
\label{eqn:(n-1)power}
    \Omega_u^{n-1} = \Omega_0^{n-1} + \partial \partial_J u \wedge \Omega^{n-2}.
\end{equation}

In this paper, we consider equation \eqref{eqn:func-type} under the assumption that $\Omega$ is  hyperK\"{a}hler.
Parallel to the complex case in \cite{tosatti2017monge}, equation \eqref{eqn:func-type} can be restated as equation \eqref{eqn:psh} in terms of quaternoinic Monge-Amp\`{e}re equation for $(n-1)$\nobreakdash-quaternoinic plurisubharmonic functions. Our main result is as follows.

\begin{theorem}
\label{thm:main}
    Let $(M, I, J, K, g, \Omega)$ be a compact hyperK\"{a}hler manifold of quaternionic dimension $n$, and $\Omega_h$ a strictly positive $(2, 0)$\nobreakdash-form with respect to $I$. Let $f$ be a smooth function on $M$. Then there is a unique pair 
    $(u, b) \in C^{\infty}(M, \mathbb{R}) \times \mathbb{R}$, solving
    \begin{equation}
    \label{eqn:psh}
        \big(\Omega_h + \frac{1}{n-1}(S_1(\p \p_J u)\Omega - \partial \partial_J u)\big)^n
        = e^{f + b} \Omega^n
    \end{equation}
    with
    \begin{equation}
    \label{main-pos}
        \Omega_h + \frac{1}{n-1}(S_1(\p \p_J u)\Omega - \partial \partial_J u) > 0, \quad
        \sup_M u = 0.
    \end{equation}
\end{theorem}

Here $S_1(\p \p_J u)$ is defined in Section \ref{notation} and related to the Chern Laplacian (see \eqref{eqn:intro-Sm}, \eqref{eqn:C0-laplace}). Recently
on a locally flat hyperhermitian manifold, Gentili and Zhang \cite{gentili2022fully} studied a general class of fully non-linear equations  including equation \eqref{eqn:psh} and they solved the equation assuming the existence of a flat hyperK\"{a}hler metric. Here we are able to remove the assumption on flatness.
 From Theorem \ref{thm:main} we obtain 

\begin{corollary}
\label{coro:main}
    Let $(M, I, J, K, g, \Omega)$ be a compact hyperK\"{a}hler manifold of quaternionic dimension $n$ and $g_0$ a quaternionic balanced (resp., Gauduchon; resp., strongly Gauduchon) metric on $M$ with induced $(2, 0)$\nobreakdash-form $\Omega_0$. Then for a given smooth function $f'$ on $M$, there exists a unique constant $b'$ and a unique quaternionic balanced (resp., Gauduchon; resp., strongly Gauduchon) metric $\Omega_u$ satisfying \eqref{eqn:(n-1)power} and solving \eqref{eqn:func-type}.
\end{corollary}

We obtain \textit{a priori} estimates and thus employ the continuity method to prove Theorem \ref{thm:main}. In Section \ref{notation} we give definitions and notations used throughout this paper, and explain the relation between equation \eqref{eqn:func-type} and \eqref{eqn:psh}. We derive $C^0$ estimate in Section \ref{C0} by establishing a Cherrier-type inequality without using the hyperK\"{a}hler condition. We derive $C^1$ estimate in Section \ref{C1} and $C^2$ estimate in Section \ref{pC2} and Section \ref{C2}. Then the main theorem is proved in last section. 
 
\section{Preliminaries and notation}
\label{notation}

On a hypercomplex manifold $(M, I, J, K)$, the exterior differential $d$ is decomposed into $d = \p_I + \bp_I$ with respect to complex structure $I$. For simplicity we denote by $\p = \p_I$ and $\bp = \bp_I$ . Verbitsky \cite{verbitsky02} introduced the operator $\p_J$ as the quaternionic analogue of $\bp$ by
\[
\p_J = J^{-1} \circ \bp \circ J.
\]
As in \cite{dinew2021hkt} we also define
\[
\bp_J = J^{-1} \circ \p \circ J .
\] The operators $\p, \bp, \p_J$ and $\bp_J$ satisfy the following properties.

\begin{lemma}[\cite{dinew2021hkt}, Lemma 2.12]
\label{lem:pre-commute}
    For a hypercomplex manifold $(M, I, J, K)$ the following holds.
    \begin{equation}
        \begin{gathered}
            \p^2 = \bp^2 = \p_J^2 = \bp_J^2 = 0,\\
            \p \bp + \bp \p = \p_J \bp_J + \bp_J \bp 
            = \p \p_J + \p_J \p = 0,\\
            \bp \bp_J + \bp_J \bp = \p_J\bp + \bp\p_J
            = \bp_J \p + \p \bp_J =0.
        \end{gathered}
    \end{equation}
\end{lemma}

On a hyperhermitian manifold $(M, I, J , K, g)$ of quaternionic dimension $n$, let 
\[ \Omega = \omega_J + i\omega_K\]
where $\omega_J$ and $\omega_K$ are the fundamental forms of $(g, J)$ and $(g, K)$ respectively. 
We denote by $\bigwedge_{I}^{p,q}(M)$ the $(p,q)$\nobreakdash-forms with respect to $I$, which we simply call $(p, q)$\nobreakdash-forms throughout this paper.
A form $\alpha \in \bigwedge_{I}^{2k, 0} (M)$ satisfying $J\alpha = \overline{\alpha}$ is called $J$\nobreakdash-real and denoted by $\alpha \in \bigwedge_{I, \mathbb{R}}^{2k, 0} (M)$. In particular, we have $\Omega \in \bigwedge_{I, \mathbb{R}}^{2, 0} (M)$. 

\begin{definition}
    A $J$\nobreakdash-real $(2,0)$\nobreakdash-form $\alpha$ is said to be
    positive (\emph{resp.} strictly positive) if 
    $\alpha(X, \overline{X}J) \geq 0$ (\emph{resp.}  $\alpha(X, \overline{X}J)  > 0$), 
    for any non-zero $(1,0)$\nobreakdash-vector $X$.
\end{definition}

In complex case, one can simultaneously diagonalize two hermitian matrices when one of them is positive definete. Similar result holds for $J$\nobreakdash-real $(2, 0)$\nobreakdash-forms.

\begin{lemma}[\cite{sroka2020}, Lemma 3]
\label{lem:pre-diag}
    Let $\alpha$ and $\beta$ be two $J$\nobreakdash-real $(2,0)$\nobreakdash-forms on a hyperhermitian manifold $M$ of quaternionic dimension $n$, and $\alpha$ is strictly positive. Then for each $x \in M$ there exists a basis
    $e_1, \overline{e_1}J, \dots, e_n, \overline{e_n}J$ of $T^{1,0}_{I, x}(M)$ such that
    \[
    \alpha(e_i, e_j) = \beta(e_i, e_j) 
    = \alpha(e_i, \overline{e_j}J) = \beta(e_i, \overline{e_j}J) = 0 
    \text{\, for\, } i \neq j. 
    \]
\end{lemma}


Analogous to positive definite $(n-1, n-1)$\nobreakdash-form in complex case  \cite{tosatti2017monge}, we define strictly positive $(2n-2,0)$\nobreakdash-form as follows.

\begin{definition}
    A $J$\nobreakdash-real $(2n-2, 0)$\nobreakdash-form $\Phi$ is said to be
    strictly positive if 
    $\Phi \wedge \alpha \geq 0$, for any positive $(2,0)$\nobreakdash-form $\alpha$, with equality if and only if $\alpha = 0$. 
    We denote all strictly positive $J$\nobreakdash-real $(2n-2, 0)$\nobreakdash-forms by $\bigwedge_{I, \mathbb{R}}^{2n-2, 0}(M)_{>0}$.
\end{definition}

The notion of (strictly) positive forms on hypercomplex manifolds can be found in various literature \cite{alesker2006plurisubharmonic, verbitsky10pos, verbitsky09balance}, and we refer readers to \cite{verbitsky10pos} for thorough discussions. For complex case see for instance \cite{demailly1997complex}.

As in \cite{gentili2022fully}, we define the Hodge star\nobreakdash-type operator
$\ast: \bigwedge_{I}^{p,0}(M) \to \bigwedge_{I}^{2n-p, 0}(M)$ by the relation
\[
\alpha \wedge \ast \beta = \frac{1}{n!} \langle \alpha, \beta \rangle_{g} \Omega^n,\quad
\text{for \,} \alpha, \beta \in \bigwedge\nolimits_{I}^{p,0}(M).
\]
Here, the pointwise inner product $\langle \, , \, \rangle_g$ is defined by
\[
\langle \alpha , \beta \rangle_g =
\frac{1}{p!}\sum g^{\lambda_1\overline{\mu}_1} \cdots g^{\lambda_p\overline{\mu}_p} 
\alpha_{\lambda_1 \cdots \lambda_p} \overline{\beta_{\mu_1 \cdots \mu_p}}, \quad
\text{for \,} \alpha, \beta \in \bigwedge\nolimits_{I}^{p,0}(M)
\]
where any $(p,0)$\nobreakdash-form $\alpha$ is locally written as
\[ 
\alpha = \frac{1}{p!} \sum \alpha_{\lambda_1 \cdots \lambda_p} 
dz^{\lambda_1} \wedge \cdots \wedge dz^{\lambda_p}.
\]
At a point $p \in M $ we can take $I$\nobreakdash-holomorphic coordinates $(z^0, \cdots, z^{2n-1})$ such that $(g_{\lambda \overline{\mu}})$ is the identity at $p$, then we have
\[
\ast (dz^{2i} \wedge dz^{2i + 1}) = 
dz^0 \wedge dz^1 \wedge \cdots \wedge \widehat{dz^{2i}} \wedge \widehat{dz^{2i + 1}}
\wedge \cdots \wedge dz^{2n -2} \wedge dz^{2n - 1}.
\]
It is easy to show that the operator $\ast$ maps 
    $\bigwedge_{I, \mathbb{R}}^{2, 0} (M)_{>0}$ to 
    $\bigwedge_{I, \mathbb{R}}^{2n-2, 0} (M)_{>0}$
    and vice versa.

\begin{definition}
\begin{enumerate}
    \item
    For a $J$\nobreakdash-real $(2, 0)$\nobreakdash-form locally written as 
    $\alpha = \sum_{i < j} \alpha_{ij} dz^i \wedge dz^j$, 
    define the Pfaffian of $\alpha$ locally by 
    $\textnormal{Pf}(\alpha)dz^0 \wedge \cdots \wedge dz^{2n-1} = \alpha^n$.
    \item 
    The Pfaffian of a $J$\nobreakdash-real $(2n-2, 0)$\nobreakdash-form $\Phi$ is defined by
    \[\textnormal{Pf}(\Phi) = \textnormal{Pf}(\frac{1}{(n-1)!}\ast\!\Phi).\]
\end{enumerate}
    
\end{definition}
In particular, we have for any $\alpha \in \bigwedge_{I, \mathbb{R}}^{2, 0} (M) $,
\begin{equation}
\label{eqn:Pf-power}
    \text{Pf}(\alpha^{n-1}) = \text{Pf}(\alpha)^{n-1}.
\end{equation}
In fact, computing at a point and using Lemma \ref{lem:pre-diag} we can write $\alpha = \sum \lambda_i dz^{2i} \wedge dz^{2i + 1}$. 
Define $\Lambda = \lambda_0 \cdots \lambda_{n-1}$, 
$\Lambda_i = \lambda_0 \cdots \hat{\lambda_i} \cdots \lambda_{n-1}$. 
Then $\text{Pf}(\alpha) = \Lambda$. On the other hand,
\[\alpha^{n-1} = (n-1)! \sum \Lambda_i
dz^0 \wedge dz^1 \wedge \cdots \wedge \widehat{dz^{2i}} \wedge \widehat{dz^{2i + 1}}
\wedge \cdots \wedge dz^{2n -2} \wedge dz^{2n - 1}.\]
By definition we have $\text{Pf}(\alpha^{n-1}) = \Lambda^{n-1}$. Hence \eqref{eqn:Pf-power} follows.

Also, observe that for any two $J$\nobreakdash-real $(2, 0)$\nobreakdash-forms
$\chi$ and $\eta$, we have 
\begin{equation}
\label{eq:pf*}
\frac{\chi^n}{\eta^n} = \frac{\text{Pf}(\chi)}{\text{Pf}(\eta)} 
=\frac{\text{Pf}(\ast\chi)}{\text{Pf}(\ast\eta)}.
\end{equation} 

For conveninence in later computation, we introduce the following definition.

\begin{definition}
For $\chi \in \bigwedge_{I, \mathbb{R}}^{2, 0}(M)$, define
\begin{equation}
\label{eqn:intro-Sm}
     S_m(\chi) = \frac{C_n^m \chi^m \wedge \Omega^{n-m}}{\Omega^n} 
     \quad \text{for} \quad 0 \leq m \leq n.
\end{equation}
\end{definition}
In particular for $u \in C^{\infty}(M, \mathbb{R})$ we have 
\begin{equation}
\label{eqn:C0-laplace}
S_1(\partial\partial_J u ) = \frac{1}{2} \Delta_{I, g} u.
\end{equation}
In fact, choose local coordinates such that
$\Omega = \sum_{i=0}^{n-1}dz^{2i} \wedge dz^{2i+1}$. Now $\omega_I$ takes the form
\[\omega_I = \frac{i}{2} \sum_{\alpha = 0}^{2n - 1} dz^{\alpha}\wedge d\overline{z}^{\alpha} .\]
Since 
\[
J(\p \p_J u ) = -J(\p_J \p u) = - \bp J \p u = \bp J^{-1} \p u = \bp \bp_J u 
= \overline{\p \p_J u}, 
\]
we see that $\p \p_J u$ is $J$\nobreakdash-real. Then compute 
\[
\begin{split}
\p \p_J u &= \sum \p (J^{-1}\bp u) = \sum \p (u_{\bar j} J^{-1}d\Bar{z}^j) \\
&= \sum u_{\bar j i}dz^i \wedge J^{-1} d\Bar{z}^j + \sum u_{\bar j} \p(J^{-1}d\Bar{z}^j)
= \sum u_{\bar j i} dz^i \wedge J^{-1} d\Bar{z}^j.
\end{split}
\]
The last equality above is derived from
\[
0 = (\bp \bp_J + \bp_J \bp)(z^i) = \bp J^{-1} \p J(z^i) + J^{-1} \p J \bp z^i 
= J^{-1} \p J d \bar z^i.
\]
Hence
\[S_1(\partial \partial_J u) 
= \frac{n \partial\partial_J u  \wedge \Omega^{n-1}}{\Omega^n}
= \sum_{\alpha = 0}^{2n-1} u_{\alpha \overline{\alpha}}\]
and 
\[
\Delta_{I, g} u 
= \frac{2n \partial \overline{\partial} u \wedge \omega_I^{2n -1} }{\omega_I^{2n} } 
= 2\sum_{\alpha=0}^{2n-1} u_{\alpha \overline{\alpha}} .\]
Thus equation \eqref{eqn:C0-laplace} holds.

Now let $(M, I, J, K, g, \Omega)$ be a compact hyperK\"{a}hler manifold and $\Omega_h$ a strictly positive $(2, 0)$\nobreakdash-form with respect to $I$. 
The quaternionic Monge–Amp\`{e}re equation for $(n-1)$\nobreakdash-quaternionic plurisubharmonic functions is written as
\begin{gather}
\label{eq:qfte-1}
      \big(\Omega_h + \frac{1}{n-1}(S_1(\p \p_J u)\Omega - \partial\partial_J u)\big)^n 
      = e^{f+b} \Omega^n\\
\label{eq:qfte-2}
      \Omega_h + \frac{1}{n-1}(S_1(\p \p_J u)\Omega - \partial\partial_J u) > 0, \quad
      \sup_M u = 0.
\end{gather}

For a quaternionic balanced metric (resp., Gauduchon; resp., strongly Gauduchon, for definitions and their correspondence with the complex case see \cite[Table 2]{MR3645310}) with induced $(2, 0)$\nobreakdash-form $\Omega_0$, we define $\Omega_h$ by 
\begin{equation}
\label{eqn:pre-omega_h}
    (n-1)! \ast \Omega_h = \Omega_0^{n-1} .
\end{equation}
We would like to show that a solution to equation \eqref{eq:qfte-1} gives rise to a solution to the quaternionic form type equation. 
We also need
\begin{equation}
\label{eqn:pre-ast-tOmega}
    \frac{1}{(n-1)!}\ast (\p \p_J u \wedge \Omega^{n-2}) = 
    \frac{1}{n-1}(S_1(\p \p_J u)\Omega - \p \p_J u)
\end{equation}
which can be seen by computing in local coordinates. We refer readers to \cite[p. 34]{gentili2022fully} for details. 
By \eqref{eq:pf*}, \eqref{eqn:pre-omega_h} and \eqref{eqn:pre-ast-tOmega}, we have 
\[
\begin{split}
    &\frac{\big(\Omega_h + \frac{1}{n-1}(S_1(\p \p_J u)\Omega - \partial\partial_J u)\big)^n}{\Omega^n}\\
= &\frac{\text{Pf}\big(\ast(\Omega_h + \frac{1}{n-1}(S_1(\p \p_J u)\Omega - \partial\partial_J u))\big)}{\text{Pf}(\ast\Omega)}\\
= &\frac{\text{Pf}(\Omega_0^{n-1} + \partial\partial_J u \wedge \Omega^{n-2})}{\text{Pf}(\Omega^{n-1})}.
\end{split}
\]

Now observe that a strictly positive $(2n-2, 0)$\nobreakdash-form $\Phi$ can be written as $\Phi = \phi^{n-1}$, where $\phi$ is a strictly positive $(2, 0)$ form. The proof here is almost identical to the arguments in \cite[p. 279-280]{michelsohn1982existence}.
Since  $\ast$ maps 
    $\bigwedge_{I, \mathbb{R}}^{2, 0} (M)_{>0}$ to 
    $\bigwedge_{I, \mathbb{R}}^{2n-2, 0} (M)_{>0}$, we have
\begin{equation}
\label{neq:pre-positive}
\Omega_0^{n-1} + \partial\partial_J u \wedge \Omega^{n-2} > 0.
\end{equation}
Thus there exist $\Omega_u$ such that
\[
\Omega_u^{n-1} = \Omega_0^{n-1} + \partial\partial_J u \wedge \Omega^{n-2}.
\]
Such $\Omega_u$ is quaternionic balanced (resp., Gauduchon; resp., strongly Gauduchon) and we have 
\[
e^{f + b} = 
\frac{\text{Pf}(\Omega_u^{n-1})}{\text{Pf}(\Omega^{n-1})} = 
\frac{\text{Pf}(\Omega_u)^{n-1}}{\text{Pf}(\Omega)^{n-1}} = 
\left(\frac{\Omega_u^n}{\Omega^n}\right)^{n-1}.
\]
It follows that a solution to equation \eqref{eq:qfte-1} solves
\begin{equation}
    \Omega_u^n = e^{f' + b'} \Omega^n.
\end{equation}
This proves Corollary \ref{coro:main}.

\begin{remark}
\label{rm:obata}
On a hyperhermitian manifold $(M, I, J, K, g)$, There exists a unique torsion free connection 
$\nabla^O$ called Obata connection \cite{obata1956affine} such that
\[
\nabla^O I = \nabla^O J = \nabla^O K = 0.
\]
It is well known that the hyperK\"{a}hler condition $d\Omega = 0$ is equivalent to $\nabla^O = \nabla^{LC} $ where $\nabla^{LC}$ is the Levi-Civita connection. Using Obata connection it is shown in \cite[Sect. 2.4]{dinew2021hkt} that one can choose around any point $x \in M$ local $I$\nobreakdash-holomorphic geodesic coordinates such that
the Christoffel symbol of $\nabla^O$ and the first derivatives of $J$ vanish at $x$. This property is crucial for $C^2$ estimate in Sect. \ref{C2}.
\end{remark}

From above remark, we have the following useful lemma.

\begin{lemma}
On a hyperhermitian manifold $(M, I, J, K, g)$, given a $J$\nobreakdash-real $(2,0)$\nobreakdash-form $\alpha$, for any point $x \in M$, one can choose around $x$ local $I$\nobreakdash-holomorphic coordinates such that the following relations hold at $x$, for all $i, j = 0, \dots, n-1$.
\begin{equation}
\label{eqn:pre-j-real}
\begin{split}
     &\alpha_{2i2j, p} = \overline{\alpha_{2i+12j+1, \bar p}}\, , \quad
    \alpha_{2i2j+1, p} = \overline{\alpha_{2j2i+1, \bar p}}\, ,\\
    &\alpha_{2i+12j, p} = \overline{\alpha_{2j+12i, \bar p}}\, , \quad
    \alpha_{2i+12j+1, p} = \overline{\alpha_{2i2j, \bar p}}\, ,
\end{split}
\end{equation}
if $\alpha$ is locally written as
\[
\alpha = \sum_{i<j} \alpha_{ij} dz^i \wedge dz^j, \quad \alpha_{ij} = -\alpha_{ji},
\] 
and 
\[
\alpha_{ij,p} = \frac{\p}{\p z^p} \alpha_{ij}, \quad 
\alpha_{ij, \bar p} = \frac{\p}{\p \Bar{z}^p} \alpha_{ij}.
\]
\end{lemma}
\begin{proof}
Choose local $I$\nobreakdash-holomorphic coordinates around $x$ such that at $x$, the first derivatives of $J$ vanish and
\[
Jdz^{2i} = -d\Bar{z}^{2i+1}, \quad 
Jdz^{2i+1} = d\Bar{z}^{2i}.
\]
The $J$ action on $1$\nobreakdash-forms is given by 
\[
J dz^i = J^i_{\bar k} d \bar{z}^k.
\]
Hence
\[
\begin{split}
J\alpha &= \sum_{k < l} \alpha_{kl} Jdz^k \wedge Jdz^l
= \sum_{i, j} \sum_{k<l} \alpha_{kl} J^k_{\bar i} J^l_{\bar j} d\Bar{z}^i \wedge d\Bar{z}^j \\
&= \sum_{i < j} \sum_{k<l} \alpha_{kl}(J^k_{\bar i} J^l_{\bar j} - J^k_{\bar j} J^l_{\bar i})
d\Bar{z}^i \wedge d\Bar{z}^j.
\end{split}
\]
Since the derivatives of $J$ vanish at $x$, taking $\p$ and evaluating at $x$ gives
\begin{equation}
    \label{eqn:pre-bp-Ja}
\p J \alpha = \sum_p \sum_{i < j} \sum_{k<l} 
\alpha_{kl, p}(J^k_{\bar i} J^l_{\bar j} - J^k_{\bar j} J^l_{\bar i})
dz^p \wedge d\Bar{z}^i \wedge d\Bar{z}^j.
\end{equation}
On the other hand
\begin{equation}
\label{eqn:pre-bp-ba}
\p \Bar{\alpha} = \sum_p \sum_{i<j} \overline{\alpha_{ij,}}{}_{p} 
dz^p \wedge d\Bar{z}^i \wedge d\Bar{z}^j.
\end{equation}
Notice at the point $x$
\[
J^{2i}_{\overline{2i+1}} = -1, \quad
J^{2i+1}_{\overline{2i}} = 1,
\]
and all the other $J^k_{\Bar{l}}$ vanish. Since $J\alpha = \Bar{\alpha}$, comparing components of \eqref{eqn:pre-bp-Ja} and \eqref{eqn:pre-bp-ba} we get for example, 
when $2i+1 < 2j +1$,
\[
\begin{split}
    \overline{\alpha_{2i+1 2j+1,}}{}_{p}  &= 
\sum_{k<l} \alpha_{kl, p}
(J^k_{\overline{2i+1}}J^l_{\overline{2j+1}} - J^k_{\overline{2j+1}}J^l_{\overline{2i+1}}) \\
&= \alpha_{2i 2j, p}J^{2i}_{\overline{2i+1}}J^{2j}_{\overline{2j+1}} 
= \alpha_{2i2j, p}\,.
\end{split}
\]
And when $2j+1 < 2i$,
\[
\begin{split}
    \overline{\alpha_{2j+1 2i,}}{}_{p}  &= 
\sum_{k<l} \alpha_{kl, p}
(J^k_{\overline{2j+1}}J^l_{\overline{2i}} - J^k_{\overline{2i}}J^l_{\overline{2j+1}}) \\
&= \alpha_{2j 2i+1, p}J^{2j}_{\overline{2j+1}}J^{2i+1}_{\overline{2i}} 
= - \alpha_{2j2i+1, p} = \alpha_{2i+12j, p} \,.
\end{split}
\]
By taking all the other combinations of $i, j$ we obtain \eqref{eqn:pre-j-real}.
\end{proof}

\begin{remark}
\label{rm:normal}
    Combining Lemma \ref{lem:pre-diag} and Remark \ref{rm:obata}, on a hyperhermitian manifold $(M, I, J, K, g, \Omega)$ of quaternionic dimension $n$, we can find local $I$\nobreakdash-holomorphic geodesic coordinates suth that $\Omega$ and another $J$\nobreakdash-real $(2,0)$\nobreakdash-form $\tOmega$ are simultaneously diagonalizable at a point $x \in M$, i.e. 
    \[
    \Omega = \sum_{i=0}^{n-1} dz^{2i} \wedge dz^{2i+1}, \quad
    \tOmega = \sum_{i=0}^{n-1} \tOmega_{2i2i+1} dz^{2i} \wedge dz^{2i+1},
    \]
    and the Christoffel symbol of $\nabla^O$ and first derivatives of $J$ vanish at $x$, i.e.
    \[
    J^l_{\bar k,i} =J^{\bar l}_{k,i} =J^{\bar l}_{k,\bar i}=J^l_{\bar k, \bar i} =0.
    \]
    We call such local coordinates the normal coordinates around $x$.
\end{remark}
\section{\texorpdfstring{$C^0$}{C0} Estimate}
\label{C0}

Recently Sroka \cite{sroka22sharp} obtained a sharp $C^0$ estimate for a class of PDEs given by the operator dominating the quaternionic Monge-Amp\`{e}re operator. Here we adopt a different approach for our purpose by establishing a Cherrier-type inequality and the lemmas in \cite{tosatti2017monge}. We obtain

\begin{theorem}
\label{thm:C0}
Let $(M, I, J, K, g, \Omega)$ be a compact hyperhermitian manifold of quaternionic dimension $n$, and $\Omega_h$ a strictly positive $(2, 0)$\nobreakdash-form with respect to $I$. Let $f$ be a smooth function on $M$.
    If $u$ is a solution to equation \eqref{eqn:psh} satisfying \eqref{main-pos}. Then there exists a constant $C$
    depending only on the fixed data $(I, J, K, g, \Omega, \Omega_h)$ and $f$ such that
    \[
    \sup_M |u| \leq C.
    \]
\end{theorem}

Notice that by maximum principal the constant $b$ in equation \eqref{eqn:psh} is uniformly bounded by $\sup_M |f|$, $\Omega$ and $\Omega_h$. In fact, at the maximum point of $u$, 
\[
S_1(\p \p_J u) \Omega - \p \p_Ju \leq 0.
\]
Hence by equation \eqref{eqn:psh} $b$ is bounded above. Similarly $b$ is also bounded below. Thus for simplicity we denote $f+b$ still as $f$ when doing estimates.

For convenience we denote 
\begin{equation}
\label{equ:defOmegaTilda}
\tOmega = \Omega_h + \frac{1}{n-1}(S_1(\partial\partial_J u)\Omega - \partial\partial_J u).
\end{equation}
The next lemma we need is straightforward.
\begin{lemma}
\begin{gather}
    \label{eq:S1ddju}
        S_1(\partial\partial_J u) = S_1(\tOmega) - S_1(\Omega_h)\\
     \label{eq:ddju}
         \partial\partial_J u = (n-1)\Omega_h - S_1(\Omega_h)\Omega + S_1(\tOmega)\Omega - (n-1)\tOmega.
\end{gather}
\end{lemma}
\begin{proof}
    From \eqref{equ:defOmegaTilda} we have 
    \[
        n \tOmega \wedge \Omega^{n-1} = 
        n \Omega_h \wedge \Omega^{n-1} + 
        \frac{n}{n-1}
        (S_1(\partial\partial_J u)\Omega^n - \partial\partial_J u \wedge \Omega^{n-1}).
    \]
    Namely,
    \[ 
        S_1(\tOmega) = S_1(\Omega_h) + 
        \frac{1}{n-1}(nS_1(\partial\partial_J u) - S_1(\partial\partial_J u))
        = S_1(\Omega_h) + S_1(\partial\partial_J u)
    \]
    This proves \eqref{eq:S1ddju}, and \eqref{eq:ddju} follows by inserting \eqref{eq:S1ddju} into 
    \eqref{equ:defOmegaTilda}.
\end{proof}

Define $\Omega_0$ by $(n-1)! \Omega_h =\ast \Omega_0^{n-1}$, we have the follwing
\begin{lemma}
\label{lem:C0-neq}
    There exists a uniform constant $C$ such that 
    \begin{equation}
    \label{ineq:c0}
         \partial\partial_J u \wedge 
        (2 \Omega_0^{n-1} + \partial\partial_J u \wedge \Omega^{n-2}) 
        \leq C\Omega^n
    \end{equation}
\end{lemma} 
\begin{proof}
    Using \eqref{eq:ddju} we compute
    \[
    \begin{split}
     &\:\partial\partial_J u \wedge 
    (2 \Omega_0^{n-1} + \partial\partial_J u \wedge \Omega^{n-1}) 
    \\
    = &\:2\big((n-1)\Omega_h - S_1(\Omega_h)\Omega\big)\wedge \Omega_0^{n-1} 
       - 2\big((n-1)\tOmega - S_1(\tOmega)\Omega\big)\wedge \Omega_0^{n-1} 
       \\
      &\:+ \big((n-1)\Omega_h - S_1(\Omega_h)\Omega - ((n-1)\tOmega - S_1(\tOmega)\Omega)\big)^2 \wedge \Omega^{n-2}
    \\
    = &\:2\big((n-1)\Omega_h - S_1(\Omega_h)\Omega\big)\wedge \Omega_0^{n-1}  
       +  \big((n-1)\Omega_h - S_1(\Omega_h)\Omega\big)^2 \wedge \Omega^{n-2} 
       \\
      &\:- 2\big((n-1)\tOmega - S_1(\tOmega)\Omega\big)\wedge \Omega_0^{n-1} 
         - 2\big((n-1)\Omega_h - S_1(\Omega_h)\Omega\big)
          \big((n-1)\tOmega - S_1(\tOmega)\Omega\big) \wedge \Omega^{n-2}
       \\
      &\: + (n-1)^2\tOmega^2\wedge \Omega^{n-2}  
          - 2(n-1)S_1(\tOmega)\tOmega \wedge \Omega^{n-1} 
          + S_1^2(\tOmega)\Omega^n
    \\
    \leq &\: C\Omega^n - 2(n-1)\tOmega\wedge\Omega_0^{n-1} 
             + 2S_1(\tOmega) \Omega \wedge \Omega_0^{n-1} 
             - 2(n-1)^2 \Omega_h \wedge \tOmega \wedge \Omega^{n-2}
         \\
         &\: + 2(n-1)S_1(\tOmega)\Omega_h \wedge \Omega^{n-1} 
             + 2(n-1)S_1(\Omega_h)\tOmega \wedge \Omega^{n-1}
             - 2S_1(\Omega_h)S_1(\tOmega)\Omega^n
         \\
         &\: + (n-1)^2\tOmega^2 \wedge \Omega^{n-2} 
             - 2(n-1)S_1(\tOmega)\tOmega \wedge \Omega^{n-1}
             + S_1^2(\tOmega)\Omega^n.
    \end{split}
    \]
    By definition of $S_1(\tOmega)$ and $S_{n-1}(\Omega_0)$, we have
    \begin{equation}
    \label{ineq:C0-many}
    \begin{split}
        &\:\partial\partial_J u \wedge 
        (2 \Omega_0^{n-1} + \partial\partial_J u \wedge \Omega^{n-1}) 
        \\
    \leq &\: C\Omega^n - 2(n-1)\tOmega\wedge\Omega_0^{n-1} 
          + \frac{2}{n}S_1(\tOmega)S_{n-1}(\Omega_0)\Omega^n
          - 2(n-1)^2\Omega_h\wedge\tOmega\wedge\Omega^{n-2} 
          \\
      &\: + \frac{2(n-1)}{n}S_1(\tOmega)S_1(\Omega_h)\Omega^n
          + \frac{2(n-1)}{n}S_1(\Omega_h)S_1(\tOmega)\Omega^n
          - 2S_1(\Omega_h)S_1(\tOmega)\Omega^n 
          \\
      &\: + \frac{2(n-1)}{n}S_2(\tOmega)\Omega^n
          - \frac{2(n-1)}{n}S_1^2(\tOmega)\Omega^n
          + S_1^2(\tOmega)\Omega^n
    \\
    = &\: C\Omega^n - 2(n-1)\tOmega\wedge\Omega_0^{n-1} 
          + \frac{2}{n}S_1(\tOmega)S_{n-1}(\Omega_0)\Omega^n
          - 2(n-1)^2\Omega_h\wedge\tOmega\wedge\Omega^{n-2} 
          \\
      &\: + \frac{2(n-2)}{n}S_1(\Omega_h)S_1(\tOmega)\Omega^n
          + \frac{2(n-1)}{n}S_2(\tOmega)\Omega^n 
          + \frac{2-n}{n}S_1^2(\tOmega)\Omega^n.
    \end{split}
    \end{equation}
    
    Choose local $I$\nobreakdash-holomorphic coordinates such that at a point, $\Omega = \sum_{i=0}^{n-1} dz^{2i} \wedge dz^{2i + 1}$ and $\Omega_0 = \sum_{i=0}^{n-1} \lambda_i  dz^{2i} \wedge dz^{2i + 1}$ with $\lambda_i > 0$.
    Since
    \[\Omega_h = \frac{1}{(n-1)!}\ast\Omega_0^{n-1} 
    = \sum_{i=0}^{n-1} \Lambda_i dz^{2i} \wedge dz^{2i+1} \] 
    where $\Lambda_i = \lambda_0 \cdots \hat{\lambda_i} \cdots \lambda_{n-1}$, we have $S_1(\Omega_h) = S_{n-1}(\Omega_0) = \sum_{i = 0}^{n-1}\Lambda_i$. Therefore
    \begin{equation}
    \label{eqn:C0-cancel-0}
        \frac{2}{n}S_1(\tOmega)S_{n-1}(\Omega_0)
        + \frac{2(n-2)}{n}S_1(\tOmega)S_1(\Omega_h)
        = \frac{2(n-1)}{n}S_1(\tOmega)S_{n-1}(\Omega_0).
    \end{equation}
    Now compute 
    \[
    \begin{split}
        2(n-1)\tOmega\wedge\Omega_0^{n-1} 
        &= 2(n-1)(n-1)! \sum_{i=0}^{n-1}\tOmega_{2i, 2i+1}\Lambda_i 
            dz^0 \wedge \cdots \wedge dz^{2n-1}
        \\
        2(n-1)^2\Omega_h \wedge \tOmega \wedge \Omega^{n-2}
        &= 2(n-1)(n-1)!\sum_{i=0}^{n-1}\Lambda_i (S_1(\tOmega) - \tOmega_{2i, 2i+1}) dz^0 \wedge \cdots \wedge dz^{2n-1}.
    \end{split}
    \]
    Thus
    \begin{equation}
    \label{eqn:C0-cancel-1}
        2(n-1)\tOmega\wedge\Omega_0^{n-1} 
        + 2(n-1)^2\Omega_h \wedge \tOmega \wedge \Omega^{n-2}
        = \frac{2(n-1)}{n}S_1(\tOmega)S_{n-1}(\Omega_0)\Omega^n.
    \end{equation}
    Combining \eqref{ineq:C0-many}, \eqref{eqn:C0-cancel-0} and \eqref{eqn:C0-cancel-1} we get 
    \begin{equation}
        \partial\partial_J u \wedge 
        (2 \Omega_0^{n-1} + \partial\partial_J u \wedge \Omega^{n-2}) 
        \leq C\Omega^n  + \frac{2(n-1)}{n}S_2(\tOmega)\Omega^n 
             + \frac{2-n}{n}S_1^2(\tOmega)\Omega^n.
    \end{equation}
    
    It remains to prove that the sum of the last two terms has a upper bound. The proof is analogous to that in \cite{tosatti2017monge}, which we give here for completeness. Choose local coordinates such that at a point,
    \[
    \begin{split}
        \Omega &= \sum_{i=0}^{n-1} dz^{2i} \wedge dz^{2i + 1} \\
        \tOmega &= \sum_{i=0}^{n-1} \mu_i  dz^{2i} \wedge dz^{2i + 1} \text{ with } 
        0 < \mu_0 \leq \cdots \leq \mu_{n-1} .
    \end{split}
    \]
    Then we have
    \[
    \begin{split}
        &\:2(n-1)S_2(\tOmega) + (2-n)S_1^2(\tOmega) \\
    =   &\:2(n-1)\sum_{i<j} \mu_i \mu_j - (n-2) (\sum_{i=0}^{n-1} \mu_i)^2 \\
    =   &\:-(n-2)\sum_{i=0}^{n-1}\mu_i^2 -2(n-2)\sum_{i<j}\mu_i\mu_j + 2(n-1)\sum_{i<j}\mu_i\mu_j\\
    =   &\:-(n-2)\sum_{i=1}^{n-1}\mu_i^2 + 2\sum_{1 \leq i < j \leq n-1} \mu_i\mu_j 
           -(n-2)\mu_0^2 + 2\mu_0\sum_{i=1}^{n-1}\mu_i\\
    \leq &\:-\sum_{1 \leq i < j \leq n-1}(\mu_i - \mu_j)^2 + 2\mu_0\sum_{i=1}^{n-1}\mu_i.
    \end{split}
    \]
    We want to show this quantity has a upper bound using the equation
    \[ 
        \mu_0 \cdots \mu_{n-1} = e^{f}.
    \]
    When $\mu_1 < \mu_{n-1}/2$, we have 
    $(\mu_1 - \mu_{n-1})^2 \geq \frac{1}{4}\mu_{n-1}^2$. Thus
    \[
        -\sum_{1 \leq i < j \leq n -1}(\mu_i - \mu_j)^2 + 2\mu_0\sum_{i=1}^{n-1}\mu_i
        \leq -\frac{1}{4}\mu_{n-1}^2 + C\mu_{n-1} \leq C',
    \]
    and the first inequality above is because $\mu_0$ has a uniform upper bound, being the smalest eigenvalue.
    When $\mu_1 \geq \mu_{n-1}/2$, then we have $\mu_i \geq \mu_{n-1}/2 $ for $i = 1, \cdots, n-1$.
    Hence 
    \[
        \mu_0 \leq \frac{C}{\mu_1 \cdots \mu_{n-1}} \leq \frac{C2^{n-2}}{\mu_{n-1}^{n-1}}.
    \]
    And in this case
    \[
        -\sum_{1 \leq i < j \leq n -1}(\mu_i - \mu_j)^2 + 2\mu_0\sum_{i=1}^{n-1}\mu_i
        \leq \frac{C'}{\mu_{n-1}^{n-1}}\mu_{n-1} = \frac{C'}{\mu_{n-1}^{n-2}} \leq C'.
    \]
    This proves the lemma.
\end{proof}

We now establish the Cherrier-type inequality:

\begin{lemma}
\label{lem:C0-cherrier}
    There exist uniform constants $C$ and $p_0$ such that for all $p \geq p_0$,
    \begin{equation}
    \label{neq:C0-cherrier}
        \int_M |\p e^{-\frac{pu}{2}}|^2_g \Omega^n \wedge \bOmega^n 
        \leq C p \int_M e^{-pu} \Omega^n \wedge \bOmega^n.
    \end{equation}
\end{lemma}
\begin{proof}
    By Lemma \ref{lem:C0-neq} we have
    \[
        \mathcal{I}:= \int_M e^{-pu} \partial\partial_J u \wedge 
        (2 \Omega_0^{n-1} + \partial\partial_J u \wedge \Omega^{n-2}) \wedge \bOmega^n
        \leq C \int_M e^{-pu} \Omega^n \wedge \bOmega^n.
    \]
    Interating by parts, we have
    \begin{align}
    \mathcal{I}=&-\int_{M}\partial e^{-pu}\wedge \partial_J u\wedge\Big(2\Omega_0^{n-1}+\partial\partial_J u\wedge\Omega^{n-2}\Big)\wedge {\overline\Omega}^n\notag\\
    &+\int_{M} e^{-pu} \partial_J u\wedge\partial\Big((2\Omega_0^{n-1}+\partial\partial_J u\wedge\Omega^{n-2})\wedge {\overline\Omega}^n\Big)\notag\\
    =&p\int_{M} e^{-pu}\partial u\wedge \partial_J u\wedge\Big(2\Omega_0^{n-1}+\partial\partial_J u\wedge\Omega^{n-2}\Big)\wedge {\overline\Omega}^n\notag\\
    &+\int_{M} e^{-pu} \partial_J u\wedge\Big((2\partial\Omega_0^{n-1}+\partial\partial_J u\wedge\partial\Omega^{n-2})\wedge {\overline\Omega}^n+(2\Omega_0^{n-1}+\partial\partial_J u\wedge\Omega^{n-2})\wedge \partial{\overline\Omega}^n\notag\Big)\notag\\
    =&\mathcal{I}_1+\mathcal{I}_2\notag
    \end{align}
    Since $\Omega_0^{n-1}+\partial\partial_J u\wedge\Omega^{n-2}>0$ (see \eqref{neq:pre-positive}), we obtain
    \begin{align}
    \mathcal{I}_1\ge p\int_{M} e^{-pu}\partial u\wedge \partial_J u\wedge\Omega_0^{n-1}\wedge {\overline\Omega}^n\notag
    \ge c_0p\int_{M} e^{-pu}\partial u\wedge \partial_J  u\wedge\Omega^{n-1}\wedge {\overline\Omega}^n,
    \end{align}
    where we use $\Omega_0\ge c_0^{\frac{1}{n-1}} \Omega$ for a positive constant $c_0$.

    Next we estimate $\mathcal{I}_2$. Indeed, we have
\begin{align*}
\mathcal{I}_2=&-\frac{1}{p}\int_M{\partial_J e^{-pu}\wedge\Big((2\partial\Omega_0^{n-1}+\partial\partial_J u\wedge\partial\Omega^{n-2})\wedge {\overline\Omega}^n
+(2\Omega_0^{n-1}+\partial\partial_J u\wedge\Omega^{n-2})\wedge \partial{\overline\Omega}^n\Big)}\notag\\
=&\frac{1}{p}\int_M{ e^{-pu} \Big((2\partial_J\partial\Omega_0^{n-1}+\partial\partial_J u\wedge\partial_J\partial\Omega^{n-2})\wedge {\overline\Omega}^n-(2\partial\Omega_0^{n-1}+\partial\partial_J u\wedge\partial\Omega^{n-2})\wedge \partial_J{\overline\Omega}^n\Big)}\notag\\
&+\frac{1}{p}\int_M{ e^{-pu} \Big((2\partial_J\Omega_0^{n-1}+\partial\partial_J u\wedge\partial_J\Omega^{n-2})\wedge \partial{\overline\Omega}^n+  ( 2\Omega_0^{n-1}+\partial\partial_J u\wedge\Omega^{n-2})\wedge {\partial_J\partial\overline\Omega}^n\Big) }\notag\\
=&\frac{1}{p}\int_M{ e^{-pu} \partial\partial_J u\wedge
\Big(\partial_J\partial\Omega^{n-2}\wedge{\overline\Omega}^n-\partial\Omega^{n-2}\wedge \partial_J{\overline\Omega}^n+\partial_J\Omega^{n-2}\wedge \partial {\overline\Omega}^n+\Omega^{n-2}\wedge\partial_J\partial {\overline \Omega}^n\Big)}\notag\\
&+\frac{1}{p}\int_M{ e^{-pu} \Big(2\partial_J\partial\Omega_0^{n-1}\wedge{\overline\Omega}^n-2\partial\Omega_0^{n-1}\wedge \partial_J{\overline\Omega}^n+2\partial_J\Omega_0^{n-1}\wedge \partial {\overline\Omega}^n+2\Omega_0^{n-1}\wedge\partial_J\partial {\overline \Omega}^n\Big)}\notag\\
=&\mathcal{I}_{21}+\mathcal{I}_{22}.
\end{align*}
$\mathcal{I}_{22}$ has the following estimate:
$$\mathcal{I}_{22}\ge -Cp^{-1}\int_M{e^{-pu} \Omega^n\wedge\overline{\Omega}^n}.$$
Integrating by parts, we have
\begin{align*}
\mathcal{I}_{21}=&\int_M{ e^{-pu} \partial u\wedge \partial_J u\wedge
\Big(\partial_J\partial\Omega^{n-2}\wedge{\overline\Omega}^n-\partial\Omega^{n-2}\wedge \partial_J{\overline\Omega}^n+\partial_J\Omega^{n-2}\wedge \partial {\overline\Omega}^n+\Omega^{n-2}\wedge\partial_J\partial {\overline \Omega}^n\Big)}\notag\\
+&\frac{1}{p}\int_M{ e^{-pu} \partial_J u\wedge
\Big(\partial_J\partial\Omega^{n-2}\wedge\partial{\overline\Omega}^n+\partial\Omega^{n-2}\wedge \partial\partial_J{\overline\Omega}^n+\partial\partial_J\Omega^{n-2}\wedge \partial {\overline\Omega}^n+\partial\Omega^{n-2}\wedge\partial_J\partial {\overline \Omega}^n\Big)}\notag\\
\ge& -C\int_M{ e^{-pu} \partial u\wedge \partial_J u\wedge\Omega^{n-1}\wedge\overline\Omega^n}.
\end{align*}

Therefore we obtain
\begin{align*}
\mathcal{I} \ge (c_0p-C)\int_M{ e^{-pu} \partial u\wedge \partial_J u\wedge\Omega^{n-1}\wedge\overline\Omega^n}-
\frac{C}{p}\int_M{ e^{-pu} \Omega^{n}\wedge\overline\Omega^n}\notag\\
\ge \frac{c_0p}{2}\int_M{ e^{-pu} \partial u\wedge \partial_J u\wedge\Omega^{n-1}\wedge\overline\Omega^n}-
\frac{C}{p}\int_M{ e^{-pu} \Omega^{n}\wedge\overline\Omega^n}.
\end{align*}
    Take $p_0 = (2c_0)^{-1}C$, then for all $p \geq p_0$,
    \[
    \frac{1}{p} \int_M \p e^{-\frac{pu}{2}} \wedge \p_J e^{-\frac{pu}{2}} \wedge \Omega^{n-1} \wedge \bOmega^n
    \leq C\int_M e^{-pu} \Omega \wedge \bOmega^n.
    \]
    This proves the lemma.
\end{proof}

\begin{proof}[Proof of Theorem \ref{thm:C0}]
    From Lemma \ref{lem:C0-cherrier}, we can prove the $C^0$ estimate using similar arguments as that in \cite{FuYau08} and \cite{TW-AJM, TW10JAMS, TW19} by regarding $M$ as a Hermitian manifold $(M, I, g)$.
    For completeness, we sketch the proof here.
    
    By \cite{TW10JAMS}, the Cherrier-type inequality \eqref{neq:C0-cherrier} implies
    \[
    e^{-p_0 \inf\limits_M u} \leq C \int_M e^{-p_0 u} \omega_I^{2n}.
    \]
    Then by \cite{FuYau08} or \cite{TW-AJM} there exist uniform constants $C_1$ and $\delta > 0$ such that
    \[
    |\{ u \leq \inf_M u + C_1 \}|_{\omega_I} \geq \delta.
    \]
    
    On the other hand, from $\sup_M u =0$ 
    and $ \Delta_{\omega_I} u = 2 S_1(\p \p_J u) \geq -2S_1(\Omega_h) $ (see \eqref{eq:S1ddju}), one can show that (see \cite{TW19})
    \[
    \int_M (-u)\omega_I^{2n} \leq C_2.
    \]
    
    Then we have
    \[
    - \delta \inf_M u \leq \int_{\{ u \leq \inf\limits_M u + C_1 \}} (-u + C_1) \leq C.
    \]
    This finishes the proof.
\end{proof}
\section{\texorpdfstring{$C^1$}{C1} Estimate}
\label{C1}

\begin{theorem}
    Let $u$ be a solution as in Theorem \ref{thm:main}. Then there exists a constant $C$
    depending only on the fixed data $(I, J, K, g, \Omega, \Omega_h)$ and $f$ such that
    \begin{equation}
    \label{neq:C1}
        |du|_g \leq C.
    \end{equation}
\end{theorem}

\begin{proof}
    A simple computation in local coordinates shows that 
    \[
    n \partial u \wedge \partial_J u \wedge \Omega^{n-1} = \frac{1}{4} |d u|_g^2 \Omega^n.
    \]
    Define
    \[
    \beta \coloneqq \frac{1}{4} |d u|_g^2.
    \]
    Following \cite{blocki2009gradient}, we consider 
    \[
    G = \log \beta - \varphi \circ u
    \]
    where $\varphi$ is a function to be determined. Suppose $G$ attain its maximum at $p$, and from now on we compute at the point $p$ using the normal coordinates around $p$ (see Remark \ref{rm:normal}).
    \[
    \begin{split}
        \partial G &= \frac{\partial \beta}{\beta} - \varphi' \partial u = 0; \\
        \partial_J G &= \frac{\partial_J \beta}{\beta} - \varphi' \partial_J u = 0;\\
        \partial\partial_J G &= \frac{\partial \partial_J \beta}{\beta} -
                                \frac{\partial \beta \wedge \partial_J \beta}{\beta^2} -
                                \varphi'' \partial u \wedge \partial_J u -
                                \varphi' \partial \partial_J u \\
                            &= \frac{\partial \partial_J \beta}{\beta} -
                               ((\varphi')^2 + \varphi'') \partial u \wedge \partial_J u
                               - \varphi' \partial \partial_J u.
    \end{split}
    \]
    Let 
    \begin{equation}
    \label{eqn:C1-A}
        A = S_{n-1}(\tOmega)\Omega^{n-1} - \tOmega^{n-1},
    \end{equation}
    where $\tOmega$ is as in the last section. Computing in normal coordinates shows
    \[
        A = (n-1)! \sum_{i = 0}^{n-1} (\sum_{j \neq i} 
        \frac{\tOmega_{01} \cdots  \tOmega_{2n-2\:2n-1}}{\tOmega_{2j2j+1}})
        dz^0 \wedge dz^1 \wedge \cdots \wedge \widehat{dz^{2i}} \wedge \widehat{dz^{2i + 1}}
        \wedge \cdots \wedge dz^{2n -2} \wedge dz^{2n - 1}.
    \]
    Thus $A$ is positive, and we have at point $p$ 
    \begin{equation}
    \label{eqn:C1less0}
    \begin{split}
     0 &\geq \frac{\partial \partial_J G \wedge A \wedge \bOmega^n}{\tOmega{\mathstrut}^n \wedge \bOmega^n} \\
     &= \frac{\partial \partial_J \beta \wedge A \wedge \bOmega^n}{\beta \tOmega{\mathstrut}^n \wedge \bOmega^n} - 
     ((\varphi')^2 + \varphi'') \frac{\partial u \wedge \partial_J u \wedge A \wedge \bOmega^n}{\tOmega{\mathstrut}^n \wedge \bOmega^n} -
     \varphi' \frac{\partial \partial_J u \wedge A \wedge \bOmega^n}{\tOmega{\mathstrut}^n \wedge \bOmega^n}.
    \end{split}
    \end{equation}
    
    We need to compute $\partial \partial_J \beta$. By definition of $\beta$ we have
    \[
    \beta \bOmega^n = n \overline{\partial} u \wedge \overline{\partial_J} u \wedge \bOmega^{n-1}.
    \]
    Taking $\partial_J$ of both sides and noticing $\partial_J \Omega = 0$ since $\Omega$ is hyperK\"{a}hler, we get
    \[
    \partial_J \beta \wedge \bOmega^n = n \partial_J \overline{\partial} u \wedge \overline{\partial_J} u \wedge \bOmega^{n-1} -
    n \overline{\partial} u \wedge \partial_J \overline{\partial_J} u \wedge \bOmega^{n-1}.
    \]
    Then taking $\partial$ of both sides, we get
    \[
    \begin{split}
    \partial \partial_J \beta \wedge \bOmega^n = 
    &n \partial \partial_J \overline{\partial} u \wedge \overline{\partial_J} u \wedge \bOmega^{n-1} + 
    n \partial_J \overline{\partial} u \wedge \partial \overline{\partial_J} u \wedge \bOmega^{n-1} \\
    &- n \partial \overline{\partial} u \wedge \partial_J \overline{\partial_J} u 
    \wedge \bOmega^{n-1}
    + n \overline{\partial} u \wedge \partial \partial_J \overline{\partial_J} u
    \wedge \bOmega^{n-1}.
    \end{split}
    \]
    From the equation 
    \begin{equation}
    \label{eqn:C1-main}
        \tOmega^n = e^f \Omega^n,
    \end{equation}
    by taking $\overline{\partial}$ we obtain
    \[
    A \wedge n \overline{\partial} \partial \partial_J u = 
    (n-1) (-n \tOmega^{n-1} \wedge \overline{\partial} \Omega_h + \overline{\partial} e^f \wedge \Omega^n ),
    \]
    and by taking $\overline{\partial_J}$ we obtain 
    \[
    A \wedge n \overline{\partial_J} \partial \partial_J u = 
    (n-1) (-n \tOmega^{n-1} \wedge \overline{\partial_J} \Omega_h + \overline{\partial_J} e^f \wedge \Omega^n ).
    \]
    Thus we have for the first term of \eqref{eqn:C1less0}
    \begin{equation}
    \label{eqn:C1-1}
        \partial \partial_J \beta \wedge A \wedge \bOmega^n = 
        I_1 + I_2 +
        n \partial_J \overline{\partial} u \wedge \partial \overline{\partial_J} u 
        \wedge \bOmega^{n-1} \wedge A - n \partial \overline{\partial} u 
        \wedge \partial_J \overline{\partial_J} u \wedge \bOmega^{n-1} \wedge A 
    \end{equation}
    where
    \[
    \begin{split}
        I_1 &= (n-1)(-n \tOmega^{n-1} \wedge \overline{\partial} \Omega_h + \overline{\partial} e^f \wedge \Omega^n ) 
        \wedge \overline{\partial_J} u \wedge \bOmega^{n-1}, \\
        I_2 &= (n-1)(n \tOmega^{n-1} \wedge \overline{\partial_J} \Omega_h - \overline{\partial_J} e^f \wedge \Omega^n ) 
        \wedge \overline{\partial} u \wedge \bOmega^{n-1}.
    \end{split}
    \]
   
    By direct computation, 
    \[
    \begin{split}
        \partial_J \overline{\partial} u &= \sum u_{\overline{ji}} J^{-1} d\overline{z^i} \wedge d\overline{z^j}; \\
        \partial \overline{\partial_J} u &= \sum u_{ij} dz^j \wedge J^{-1} dz^i; \\
        \partial \overline{\partial} u &= \sum u_{i\overline{j}} dz^i \wedge d\overline{z^j}; \\
        \partial_J \overline{\partial_J} u &= \sum u_{i\overline{j}} J^{-1} d\overline{z^j} \wedge J^{-1} dz^i;
    \end{split}
    \]
    Thus the third term of \eqref{eqn:C1-1} become
    \begin{equation}
    \label{eqn:C1-1-3}
        n \partial_J \overline{\partial} u \wedge \partial \overline{\partial_J} u 
        \wedge \bOmega^{n-1} \wedge A =
        \frac{1}{n} \sum_{k=0}^{n-1} \sum_{j=0}^{2n-1} (\sum_{i\neq k}\frac{1}{\tOmega_{2i2i+1}} )(|u_{2kj}|^2 + |u_{2k+1j}|^2)
        \tOmega^n \wedge \bOmega^n;
    \end{equation}
    and the forth term
    \begin{equation}
    \label{eqn:C1-1-4}
    -n \partial \overline{\partial} u \wedge \partial_J \overline{\partial_J} u 
        \wedge \bOmega^{n-1} \wedge A =
        \frac{1}{n} \sum_{k=0}^{n-1} \sum_{j=0}^{2n-1} (\sum_{i\neq k}\frac{1}{\tOmega_{2i2i+1}} ) 
        (|u_{2k\overline{j}}|^2 + |u_{2k+1\overline{j}}|^2)
        \tOmega^n \wedge \bOmega^n.
    \end{equation}
    For $I_1$ and $I_2$ we have
    \begin{equation}
    \label{eqn:C1-1-1}
    \begin{split}
    \frac{1}{n-1} I_1 &= -n \tOmega^{n-1} \wedge \overline{\partial} \Omega_h \wedge \overline{\partial_J} u \wedge \bOmega^{n-1}  - \overline{\partial_J} u \wedge  \overline{\partial} e^f \wedge \Omega^n \wedge \bOmega^{n-1} \\
    &= - \frac{1}{n} \sum_{i = 0}^{n-1} \sum_{j = 0}^{2n-1} 
    \frac{(\Omega_h)_{2i2i+1, \overline{j}} u_{j}  }{\tOmega_{2i2i+1}}\tOmega^n \wedge \bOmega^n + 
    \frac{1}{n} \sum_{j=0}^{2n-1} \frac{u_j(e^f)_{\overline{j}}}{e^f} \tOmega^n \wedge \bOmega^n
    \end{split}
    \end{equation}
    and 
    \begin{equation}
    \label{eqn:C1-1-2}
    \begin{split}
        \frac{1}{n-1} I_2 &= n \tOmega^{n-1} \wedge \overline{\partial_J} \Omega_h \wedge \overline{\partial} u \wedge \bOmega^{n-1}  + \overline{\partial} u \wedge  \overline{\partial_J} e^f \wedge \Omega^n \wedge \bOmega^{n-1} \\
    &= - \frac{1}{n} \sum_{i = 0}^{n-1} \sum_{j = 0}^{2n-1} 
    \frac{(\bOmega_h)_{2i2i+1, j} u_{\overline{j}}  }{\tOmega_{2i2i+1}}\tOmega^n \wedge \bOmega^n + 
    \frac{1}{n} \sum_{j=0}^{2n-1} \frac{u_{\overline{j}}(e^f)_{j}}{e^f} \tOmega^n \wedge \bOmega^n.
    \end{split}
    \end{equation}
    Combining \eqref{eqn:C1-1-1}, \eqref{eqn:C1-1-2}, \eqref{eqn:C1-1-3}, \eqref{eqn:C1-1-4} we obtain estimate of \eqref{eqn:C1-1}
    \begin{equation}
    \label{eqn:C1-1-local}
    \begin{split}
       \frac{\p \p_J \beta \wedge A \wedge \bOmega^n}{\beta \tOmega{\mathstrut}^n \wedge \bOmega^n}
       &= - \frac{1}{n \beta} \sum_{i = 0}^{n-1} \sum_{j = 0}^{2n-1} 
    \frac{(\Omega_h)_{2i2i+1, \overline{j}} u_{j} + (\bOmega_h)_{2i2i+1, j} u_{\overline{j}} }
    {\tOmega_{2i2i+1}} 
    + \frac{1}{n \beta} \sum_{j=0}^{2n-1} 
    \frac{u_j(e^f)_{\overline{j}} + u_{\overline{j}}(e^f)_{j} }{e^f} \\
    &+ \frac{1}{n \beta} \sum_{k=0}^{n-1} \sum_{j=0}^{2n-1} \sum_{i\neq k}
    \frac{|u_{2kj}|^2 + |u_{2k+1j}|^2}{\tOmega_{2i2i+1}} 
    + \frac{1}{n \beta} \sum_{k=0}^{n-1} \sum_{j=0}^{2n-1} \sum_{i\neq k}
    \frac{|u_{2k\overline{j}}|^2 + |u_{2k+1\overline{j}}|^2}{\tOmega_{2i2i+1}} .
    \end{split} 
    \end{equation}
    Again by direct computation, the second term of \eqref{eqn:C1less0} is
    \begin{equation}
    \label{eqn:C1-2}
    \partial u \wedge \partial_J u \wedge A \wedge \bOmega^n =
     \frac{1}{n} \sum_{i=0}^{n-1}  (\sum_{k\neq i}\frac{1}{\tOmega_{2k2k+1}} ) 
     (|u_{2i}|^2 + |u_{2i+1}|^2) \tOmega^n \wedge \bOmega^n.
    \end{equation}
    For the third term of \eqref{eqn:C1less0}, we compute
    \begin{equation}
    \begin{split}
        \partial \partial_J u \wedge A &= \partial \partial_J u 
        \wedge (\frac{n \tOmega^{n-1} \wedge \Omega}{\Omega^n}\Omega^{n-1} -\tOmega^{n-1}) \\
        &= (S_1(\partial \partial_J u)\Omega - \partial \partial_J u ) \wedge \tOmega^{n-1} \\
        &= (n-1) (\tOmega^n - \Omega_h \wedge \tOmega^{n-1}).
    \end{split}
    \end{equation}
    By compactness of $M$, there exists $\epsilon >0$ such that $\Omega_h \geq \epsilon \Omega$, we obtain
    \begin{equation}
    \label{neq:C1-3}
    \begin{split}
    - \varphi' \frac{\partial \partial_J u \wedge A \wedge \bOmega^n}{\tOmega{\mathstrut}^n \wedge \bOmega^n} &= -(n-1)\varphi'  +(n-1)\varphi' \frac{\Omega_h \wedge \tOmega^{n-1} \wedge \bOmega^n }{\tOmega{\mathstrut}^n \wedge \bOmega^n} \\
    &\geq -(n-1)\varphi' + 
    \frac{\epsilon (n-1) \varphi'}{n} \sum_{i=0}^{n-1} \frac{1}{\tOmega_{2i2i+1}}.
    \end{split}
    \end{equation}
    We may assume $\beta \gg 1$ otherwise we are finished.
    The inequality \eqref{eqn:C1less0} become
    \begin{equation}
    \label{neq:C1-final}
    \begin{split}
        0 \geq & \frac{n-1}{n \beta e^f} \sum_{i=0}^{2n-1} (u_i(e^f)_{\overline{i}} + u_{\overline{i}}(e^f)_i) \\
        &- \frac{(\varphi')^2 + \varphi''}{n}\sum_{i=0}^{n-1} (\sum_{k\neq i}\frac{1}{\tOmega_{2k2k+1}} )  (|u_{2i}|^2 + |u_{2i+1}|^2) \\
        &- (n-1)\varphi' + \frac{n-1}{n}(\epsilon\varphi' - C_1 \frac{\sum u_j}{\beta} - C_2\frac{\sum u_{\overline{j}}}{\beta})\sum_{i=0}^{n-1}\frac{1}{\tOmega_{2i2i+1}}.
    \end{split}
    \end{equation}
    The first term is bounded from below. Now we take 
    \begin{equation}
        \varphi(t) = \frac{\log(2t+C_0)}{2}.
    \end{equation}
    where $C_0$ is determined by $C^0$ estimate, and rewrite \eqref{neq:C1-final} as
    \begin{equation}
    \label{neq:C1-bound}
        C_3 \geq C_4\sum_{i=0}^{n-1} (\sum_{k\neq i}\frac{1}{\tOmega_{2k2k+1}}) 
        (|u_{2i}|^2 + |u_{2i+1}|^2) + 
        C_5 \sum_{i=0}^{n-1}\frac{1}{\tOmega_{2i2i+1}}.
    \end{equation}
    Thus for any fixed $i$ 
    \[
    \tOmega_{2i2i+1} \geq \frac{C_5}{C_3} \geq C.
    \]
    By equation \eqref{eqn:C1-main} we also have
    \[
    \frac{1}{\tOmega_{2i2i+1} } = e^{-f} \prod_{j\neq i}\tOmega_{2j2j+1} \geq \frac{C^{n-1}}{\sup_{M}e^f}.
    \]
    From the bound on all $\tOmega_{2i2i+1}$, we obtain the bound on $\beta$ by \eqref{neq:C1-bound}.
\end{proof}
\section{Bound on \texorpdfstring{$\partial \partial_J u$}{∂∂ju}}
\label{pC2}

\begin{theorem}
    Let $u$ be a solution as in Theorem \ref{thm:main}. Then there exists a constant $C$
    depending only on the fixed data $(I, J, K, g, \Omega, \Omega_h)$ and $f$ such that
    \begin{equation}
    \label{neq:pC2}
        |\partial \partial_J u|_g \leq C.
    \end{equation}
\end{theorem}
\begin{proof}
    For simplicity denote
    \[
        \eta  =S_1(\p \p_J u).
    \]
    Consider the function 
    \[
        G = \log \eta - \varphi \circ u
    \]
    where the function $\varphi$ is as in the previous section. We compute at a maximum point $p$ of $G$ using the normal coordinates around $p$ (see Remark \ref{rm:normal}). We have
    \[
    \begin{split}
        \p G &= \frac{\p \eta }{\eta} - \varphi' \p u = 0; \\
        \p_J G &= \frac{\p_J \eta}{\eta} - \varphi' \p_J u = 0; \\
        \p \p_J G &= \frac{\p \p_J \eta}{\eta} - ((\varphi')^2 + \varphi'')\p u \wedge \p_J u
                    - \varphi' \p \p_J u.
    \end{split}
    \]
    Let $A$ be as before (see \eqref{eqn:C1-A}), then at point $p$ we have
    \begin{equation}
    \label{neq:pC2-beginning}
        \begin{split}
            0 &\geq \frac{\p \p_J G \wedge A \wedge \bOmega^n}{\tOmega{\mathstrut}^n \wedge \bOmega^n} \\
            &= \frac{\p \p_J \eta \wedge A \wedge \bOmega^n}{\eta \tOmega{\mathstrut}^n \wedge \bOmega^n}
            - ((\varphi')^2 + \varphi'')\frac{\p u \wedge \p_J u \wedge A \wedge \bOmega^n }{\tOmega{\mathstrut}^n \wedge \bOmega^n}
            - \varphi' \frac{\p \p_J u \wedge A \wedge \bOmega^n}{\tOmega{\mathstrut}^n \wedge \bOmega^n}.
        \end{split}
    \end{equation}
    The second and the third term were dealt with in the previous section. We now focus on $\p \p_J \eta$ in the first term. 
    
    By definition $\eta$ is real, and  
    \[
    \eta \bOmega^n = n \bp \bp_J  u \wedge \bOmega^{n-1}.
    \]
    Under the hyperK\"{a}hler condition $\dmath \Omega = 0$,  differentiating twice the above equation gives
    \begin{equation}
        \p \p_J \eta \wedge \bOmega^n 
        = n  \p \p_J \bp \bp_J  u  \wedge \bOmega^{n-1}
        = n \bp \bp_J  \p \p_J  u  \wedge \bOmega^{n-1}
    \end{equation}
    The last equality above is due to Lemma \ref{lem:pre-commute}.
    
    We know that (recall \eqref{eq:ddju})
    \[
    \partial\partial_J u = (n-1)\Omega_h - S_1(\Omega_h)\Omega + S_1(\tOmega)\Omega - (n-1)\tOmega.
    \]
    Thus 
    \begin{equation}
        \bp \bp_J \p \p_J u = (n-1)\bp \bp_J \Omega_h - 
        \bp \bp_J S_1(\Omega_h)\wedge \Omega +
        \bp \bp_J S_1(\tOmega)\wedge \Omega 
        - (n-1)\bp \bp_J \tOmega.
    \end{equation}
    Here we again used the hyperK\"{a}hler condition on $\Omega$. Now we have
    \begin{equation}
    \label{eqn:pC2-forth-diff}
    \begin{split}
        \p \p_J \eta \wedge A \wedge \bOmega^n
        &= n A \wedge  \bp \bp_J  \p \p_J u \wedge \bOmega^{n-1} \\
        &= n(n-1) A \wedge \bp \bp_J \Omega_h \wedge \bOmega^{n-1} 
        - n\bp \bp_J S_1(\Omega_h) \wedge A \wedge \Omega \wedge \bOmega^{n-1}\\
        &\quad + n\bp \bp_J S_1(\tOmega)\wedge A \wedge \Omega \wedge \bOmega^{n-1}
        - n(n-1) A \wedge \bp \bp_J \tOmega \wedge \bOmega^{n-1}
    \end{split}
    \end{equation}
    Notice that
    \[
        A \wedge \Omega = S_{n-1}(\tOmega)\Omega^n - \tOmega^{n-1} \wedge \Omega 
        = \frac{n-1}{n}S_{n-1}(\tOmega)\Omega^n
    \]
    and 
    \[
        \bp \bp_J S_1(\tOmega) \wedge \Omega^n 
        = n \bp \bp_J \tOmega \wedge \Omega^{n-1}.
    \]
    The third term of \eqref{eqn:pC2-forth-diff} becomes
    \[
    \begin{split}
        \bp \bp_J S_1(\tOmega)\wedge A \wedge \Omega \wedge \bOmega^{n-1}
        &= \bp \bp_J S_1(\tOmega)\wedge (\Omega^n 
        \cdot \frac{n-1}{n}S_{n-1}(\tOmega)) \wedge \bOmega^{n-1} \\
        &= (n-1) S_{n-1}(\tOmega) \bp \bp_J \tOmega 
        \wedge \Omega^{n-1} \wedge \bOmega^{n-1}.
    \end{split}
    \]
    The forth term is 
    \[
        A \wedge \bp \bp_J \tOmega \wedge \bOmega^{n-1}
        =  S_{n-1}(\tOmega) \bp \bp_J \tOmega 
        \wedge \Omega^{n-1} \wedge \bOmega^{n-1} 
        - \tOmega^{n-1} \wedge \bp \bp_J \tOmega \wedge \bOmega^{n-1}.
    \]
    The first two terms of \eqref{eqn:pC2-forth-diff} are similar and we get 
    \[
        \p \p_J \eta \wedge A \wedge \bOmega^n 
        = n(n-1) \bp \bp_J \tOmega \wedge  \tOmega^{n-1} \wedge \bOmega^{n-1} 
          - n(n-1) \bp \bp_J \Omega_h \wedge  \tOmega^{n-1} \wedge \bOmega^{n-1} 
    \]
    and
    \begin{equation}
    \label{eqn:pC2-target}
    \begin{split}
        \frac{\p \p_J \eta \wedge A \wedge \bOmega^n}{\eta \tOmega{\mathstrut}^n \wedge \bOmega^n}
        &= n(n-1)\frac{
            \bp \bp_J \tOmega  \wedge \tOmega^{n-1}
            \wedge \bOmega^{n-1}}{
            \eta \tOmega{\mathstrut}^n \wedge \bOmega^n
            } 
            - n(n-1)\frac{
            \bp \bp_J \Omega_h  \wedge \tOmega^{n-1}
            \wedge \bOmega^{n-1}}{
            \eta \tOmega{\mathstrut}^n \wedge \bOmega^n
            } \\
        &= \frac{n-1}{\eta n} \sum_{i=0}^{n-1} \sum_{p=0}^{2n-1}
        \frac{\tOmega_{2i2i+1,p\bar p}}{\tOmega_{2i2i+1}}
        - \frac{n-1}{\eta n} \sum_{i=0}^{n-1} \sum_{p=0}^{2n-1}
        \frac{(\Omega_h)_{2i2i+1,p\bar p}}{\tOmega_{2i2i+1}}\\
        &\geq \frac{n-1}{\eta n} \sum_{i=0}^{n-1} \sum_{p=0}^{2n-1}
        \frac{\tOmega_{2i2i+1,p\bar p}}{\tOmega_{2i2i+1}}
        - \frac{C_1}{\eta} \sum_{i=0}^{n-1} \frac{1}{\tOmega_{2i2i+1}}.
    \end{split}
    \end{equation}
    We now rewrite the right hand side of \eqref{eqn:pC2-target} using the equation
    \begin{equation}
    \label{eqn-log}
        \text{Pf} (\tOmega_{ij}) = e^f \text{Pf} (\Omega_{ij}).
    \end{equation}
    Take logarithm of both sides
    \begin{equation}
    \label{eqn:pC2-log}
        \log \text{Pf} (\tOmega_{ij}) = f + \log \text{Pf} (\Omega_{ij}).
    \end{equation}
    Since $\Omega^n = \text{Pf}(\Omega_{ij})dz^0 \wedge \cdots \wedge dz^{2n-1}$ and 
    $\bp \Omega = 0$, we have $\bp \text{Pf}(\Omega) = 0$. Taking $\bp$ of \eqref{eqn:pC2-log}, since $\text{Pf}(\tOmega_{ij})^2 = \det(\tOmega_{ij})$, we get
    \begin{equation}
        \frac{1}{2} \sum \tOmega^{ij} \tOmega_{ji,\bar p} = f_{\bar p}.
    \end{equation}
    Taking $\p$ of both sides we obtain
    \begin{equation}
    \label{eqn:pC2-2th-eq-diff}
        \frac{1}{2} \sum \tOmega^{ij} \tOmega_{ji,\bar p p}= 
        \frac{1}{2} \sum \tOmega^{ik} \tOmega _{kl,p} \tOmega^{lj} \tOmega_{ji,\bar p} + f_{p \bar p} .
    \end{equation}
    Writing in local coordinates, the left hand side of \eqref{eqn:pC2-2th-eq-diff} is
    \begin{equation}
        \label{eqn:pC2-eq-diff-left}
        \frac{1}{2} \sum \tOmega^{2i2i+1} \tOmega_{2i+12i,p \bar p} 
        + \frac{1}{2} \sum \tOmega^{2i+1 2i} \tOmega_{2i2i+1,p \bar p} 
        = \sum \frac{\tOmega_{2i2i+1, p \bar p}}{\tOmega_{2i2i+1}}.
    \end{equation}
    We claim that the first term of the right hand side of \eqref{eqn:pC2-2th-eq-diff} is positive, i.e.
    \begin{equation}
    \label{neq:pC2-claim}
        \sum \tOmega^{ik} \tOmega _{kl,p} \tOmega^{lj} \tOmega_{ji,\bar p} \geq 0.
    \end{equation}
    Indeed, in canonical coordinates,
    \[
            \begin{split}
                \sum \tOmega^{ik} \tOmega _{kl,p} \tOmega^{lj} \tOmega_{ji,\bar p} 
                &= \tOmega^{2i2i+1}(\tOmega^{2j2j+1}\tOmega_{2i+12j,p}\tOmega_{2j+12i,\bar p}
                +\tOmega^{2j+12j}\tOmega_{2i+12j+1,p}\tOmega_{2j2i,\bar p})\\
                &\quad +\tOmega^{2i+12i}
                (\tOmega^{2j2j+1}\tOmega_{2i2j,p}\tOmega_{2j+12i+1, \bar p}
                + \tOmega^{2j+12j}\tOmega_{2i2j+1,p}\tOmega_{2j2i+1, \bar p})\\
                &= \sum \frac{\tOmega_{2i+12j,p}\tOmega_{2j+12i,\bar p}
                    + \tOmega_{2i2j+1,p}\tOmega_{2j2i+1, \bar p} 
                    }{\tOmega_{2i2i+1}\tOmega_{2j2j+1}} \\
                &\quad - \sum \frac{\tOmega_{2i+12j+1,p}\tOmega_{2j2i,\bar p}
                    + \tOmega_{2i2j,p}\tOmega_{2j+12i+1, \bar p} 
                    }{\tOmega_{2i2i+1}\tOmega_{2j2j+1}}  
            \end{split}
    \]
    Since $\tOmega$ is $J$\nobreakdash-real, using relation \eqref{eqn:pre-j-real} we see that
    \begin{equation}
        \tOmega^{ik} \tOmega _{klp} \tOmega^{lj} \tOmega_{ji\bar p} 
        = \sum \frac{|\tOmega_{2i+12j,p}|^2 + |\tOmega_{2i2j+1,p}|^2
        + |\tOmega_{2i+12j+1,p}|^2 + |\tOmega_{2i2j,p}|^2
        }{\tOmega_{2i2i+1}\tOmega_{2j2j+1}}
    \end{equation}
    therefore \eqref{neq:pC2-claim} holds. By \eqref{eqn:pC2-target}, \eqref{eqn:pC2-2th-eq-diff}, \eqref{eqn:pC2-eq-diff-left} and \eqref{neq:pC2-claim} we have
    \begin{equation}
        \frac{\p \p_J \eta \wedge A \wedge \bOmega^n}{\eta \tOmega{\mathstrut}^n \wedge \bOmega^n} 
        \geq \frac{n-1}{2 \eta n } \Delta_{I, g}f
        - \frac{C_1}{\eta} \sum_{i=0}^{n-1} \frac{1}{\tOmega_{2i2i+1}}.
    \end{equation}
    By \eqref{eqn:C1-2} and \eqref{neq:C1-3}, the inequality \eqref{neq:pC2-beginning} now becomes
    \begin{equation}
    \label{neq:pC2-last}
    \begin{split}
        0 \geq &\frac{n-1}{2 \eta n} \Delta_{I, g}f 
        - \frac{(\varphi')^2 + \varphi''}{n}\sum_{i=0}^{n-1} (\sum_{k\neq i}\frac{1}{\tOmega_{2k2k+1}} )  
        (|u_{2i}|^2 + |u_{2i+1}|^2) \\
        &-(n-1)\varphi' + 
        \left( \frac{\epsilon (n-1) \varphi'}{n}- \frac{C_1}{\eta} \right)
        \sum_{i=0}^{n-1} \frac{1}{\tOmega_{2i2i+1}}.
    \end{split}
    \end{equation}
    Assuming $\eta \gg 1 $, we obtain from \eqref{neq:pC2-last}
    \begin{equation}
        C_2 \geq 
        C_3 \sum_{i=0}^{n-1}\frac{1}{\tOmega_{2i2i+1}}
    \end{equation}
    and hence all $\tOmega_{2i2i+1}$ are uniformly bounded. 
    Since $\eta = S_1(\p \p_J u) = S_1(\tOmega) - S_1(\Omega_h) $, we can therefore obtain a unform bound on $\eta$.
\end{proof}
\section{\texorpdfstring{$C^2$}{C2} Estimate}
\label{C2}

\begin{theorem}
    Let $u$ be a solution as in Theorem \ref{thm:main}. Then there exists a constant $C$
    depending only on the fixed data $(I, J, K, g, \Omega, \Omega_h)$ and $f$ such that
    \begin{equation}
    \label{neq:C2}
        |\nabla^2 u|_g \leq C.
    \end{equation}
\end{theorem}
\begin{proof}
    Since the sum of eigenvalues of $\nabla^2 u $ is bounded below by
    \[
    \frac{1}{2} \Delta_{I, g} u = S_1(\p \p_J u) = S_1(\tOmega) - S_1(\Omega_h) 
    \geq - S_1(\Omega_h),
    \]
    it is sufficient to show that the maximum eigenvalue is bounded from above.
    Define a function on $M$ as in \cite{blocki2013complex}
    \[
    \lambda(x) =  \sup_{X \in S(T_xM)} (\nabla^2 u)(X, X)
    \]
    where $S(T_xM)$ denotes unit tangent vectors at $x$.  

    Consider the function 
    \[
    G = \lambda + \frac{1}{4} |du|_g^2.
    \]
    Since we have obtained $C^1$ estimate, it is sufficient to estimate $G$ at a maximum point $p \in M$.  
    In the normal coordinates around $p$ we introduce real coordinates \begin{equation}
    \label{eqn:C2-real-coor}
        z^j = t_j + it_{2n + j}, \quad j = 0, \cdots, 2n-1,
    \end{equation}
    and compute
    \begin{equation}
    \label{eqn:C2-nabla2u}
        \nabla^2 u = \nabla(u_{t_j} d t_j) 
        = u_{t_i t_j} d t_i \otimes d t_j - \Gamma^k_{ji} u_{t_j} d t_i \otimes d t_k,
    \end{equation}
    where $\Gamma^k_{ji}$ is the Christoffel symbol of $\nabla$ with respect to 
    $\{\frac{\p}{\p t_j}\}^{4n-1}_{i=0}$.
    Suppose
    \[
    X(p) = \sum_{j = 0}^{4n-1} X^j(p) \frac{\p}{\p t_j}(p)
    \]
    is the vector realizing the supremum of $\nabla^2 u$ at $p$, and we extend it to a constant vector field $X$ near $p$, i.e.
    \[
    X = \sum_{j=0}^{4n-1} X^j(p) \frac{\p}{\p t_j}.
    \]
    Then define in a sufficiently small neighbourhood,
    \[
    \begin{split}
        \Tilde{\lambda} &= \nabla^2 u (X, X) \\
        \Tilde{G} &= \Tilde{\lambda} + \frac{1}{4}|du|_g^2.
    \end{split}
    \]
    Notice that $\Tilde{\lambda} \leq \lambda$, $\Tilde{\lambda}(p) = \lambda(p)$. 
    Hence $\Tilde{G}$ also attain its maximum at $p$ near $p$, and $\Tilde{\lambda}$ therefore $\Tilde{G}$ is smooth near $p$.
    By \eqref{eqn:C2-nabla2u} we have
    \begin{equation}
    \label{eqn:C2-tilde-lambda}
        \Tilde{\lambda} = D^2_X u - \Gamma^k_{ji} u_{t_j} X^i X^k
    \end{equation}
    where $D$ denotes the usual derivative with respect to real coordinates.
    
    Let $A$ be as before (see \eqref{eqn:C1-A}), then at the point $p$ we get 
    \begin{equation}
    \label{neq:C2-two-terms}
            0 \geq
            \frac{\p \p_J \Tilde{G} \wedge A \wedge \bOmega^n}{\tOmega{\mathstrut}^n \wedge \bOmega^n} 
            = \frac{\p \p_J \Tilde{\lambda} \wedge A \wedge \bOmega^n}{\tOmega{\mathstrut}^n \wedge \bOmega^n}
            + 
            \frac{\frac{1}{4}\p \p_J |du|_g^2\wedge A \wedge \bOmega^n}{\tOmega{\mathstrut}^n \wedge \bOmega^n}.
    \end{equation}
    In local coordinates, the first term is
    \begin{equation}
    \label{eqn:C2-first-term}
    \begin{split}
        \frac{\p \p_J \Tilde{\lambda} \wedge A \wedge \bOmega^n}{\tOmega{\mathstrut}^n \wedge \bOmega^n}
        &= \frac{1}{n} \sum_{p=0}^{n-1} \sum_{i \neq p} 
        \frac{\Tilde{\lambda}_{2p\overline{2p}} 
        + \Tilde{\lambda}_{2p+1\overline{2p+1}}}{\tOmega_{2i2i+1}} \\
        &= \frac{1}{n} \sum_{i=0}^{n-1} \sum_{p \neq i} 
        \frac{\Tilde{\lambda}_{2p\overline{2p}} 
        + \Tilde{\lambda}_{2p+1\overline{2p+1}}}{\tOmega_{2i2i+1}}.
    \end{split}
    \end{equation}

    Differentiating \eqref{eqn:C2-tilde-lambda} twice gives
    \begin{equation}
    \label{neq:C2-first-2diff}
        \begin{split}
            \Tilde{\lambda}_{p \overline{p}} &= 
            D^2_{X} u _{p \overline p} - \Gamma^k_{ji p \overline p}u_{t_j} X^i X^k - 
            \Gamma^k_{jip} u_{t_j t_{\overline p}} X^i X^k - \Gamma^k_{ji \overline p} u_{t_j t_p} X^i X^k \\
            &\geq D^2_X u_{p \overline p} - C_1(\Tilde{\lambda} + 1).
        \end{split}
    \end{equation}
    Here we used Remark \ref{rm:normal} and the fact that derivatives of $\Gamma^k_{ij}$ depend only on $g$, and the gradient of $u$ is bounded. In addition 
    \[
    |u_{t_i t_j}| \leq C_2(1 + \Tilde{\lambda}).
    \]

    By \eqref{neq:pC2} and \eqref{eq:S1ddju} we know that
    \begin{equation}
    \label{neq:C2-tOmega-bound}
        \frac{1}{C_3} \leq \tOmega_{2i2i+1} \leq C_3.
    \end{equation}

    Applying \eqref{neq:C2-first-2diff} and \eqref{neq:C2-tOmega-bound} we can estimate 
    \eqref{eqn:C2-first-term}:
    \begin{equation}
    \label{neq:C2-first-term-local}
    \begin{split}
         \frac{1}{n} \sum_{i=0}^{n-1} \sum_{p \neq i} 
        \frac{\Tilde{\lambda}_{2p\overline{2p}} 
        + \Tilde{\lambda}_{2p+1\overline{2p+1}}}{\tOmega_{2i2i+1}}
        &\geq 
        \frac{1}{n}\sum_{i=0}^{n-1} \sum_{p \neq i} 
        \frac{D^2_X u_{2p\overline{2p}} + D^2_X u_{2p+1\overline{2p+1}}}{\tOmega_{2i2i+1}} 
        -C_1(\Tilde{\lambda} +1) \\
        &\geq 
        C_4 \sum_{p = 0}^{2n-1}D^2_X u_{p \overline{p}} - C_1(\Tilde{\lambda} +1)
    \end{split}
    \end{equation}
    To deal with the first term of the right hand side, we use equation \eqref{eqn-log}
    \[
    \log \text{Pf} (\tOmega_{ij}) = f + \log \text{Pf} (\Omega_{ij}).
    \]
    Differentiating twice in direction $X$, we get 
    \begin{equation}
    \label{eqn:C2-3rd-diff}
        \frac{1}{2} \sum \tOmega^{ij} D^2_X \tOmega_{ji}= 
        \frac{1}{2} \sum \tOmega^{ik} D_X \tOmega _{kl} \tOmega^{lj} D_X \tOmega_{ji} + D^2_X f + D^2_X \log \text{Pf}(\Omega_{ij}).
    \end{equation}
    As in previous section, 
    \[
    \begin{split}
        \sum \tOmega^{ik} D_X \tOmega _{kl} \tOmega^{lj} D_X \tOmega_{ji} =
        &\sum \frac{D_X\tOmega_{2i+12j}D_X\tOmega_{2j+12i}
                    + D_X\tOmega_{2i2j+1}D_X\tOmega_{2j2i+1} 
                    }{\tOmega_{2i2i+1}\tOmega_{2j2j+1}} \\
               &- \sum \frac{D_X\tOmega_{2i+12j+1}D_X\tOmega_{2j2i}
                    + D_X\tOmega_{2i2j}D_X\tOmega_{2j+12i+1} 
                    }{\tOmega_{2i2i+1}\tOmega_{2j2j+1}}.
    \end{split}
    \]
    Notice that for $p = 0, \dots, 2n-1$,
    \[
    \frac{\p}{\p t_p}  \tOmega_{ij} 
    = \frac{\p}{\p z^p} \tOmega_{ij} + \frac{\p}{\p \Bar{z}^p} \tOmega_{ij},
    \quad 
    \frac{\p}{\p t_{2n+p}}  \tOmega_{ij} 
    = -i\big(\frac{\p}{\p \Bar{z}^p} \tOmega_{ij} - \frac{\p}{\p z^p} \tOmega_{ij}\big).
    \]
    Hence by \eqref{eqn:pre-j-real}, we obtain
    \[
    \begin{split}
     &D_X\tOmega_{2i2j} = \overline{D_X\tOmega_{2i+12j+1}}\, , \quad
    D_X\tOmega_{2i2j+1} = \overline{D_X\tOmega_{2j2i+1}}\, ,\\
    &D_X\tOmega_{2i+12j} = \overline{D_X\tOmega_{2j+12i}}\, , \quad
    D_X\tOmega_{2i+12j+1} = \overline{D_X\tOmega_{2i2j}}\, .
    \end{split}
    \]
    Therefore
    \[
    \sum \tOmega^{ik} D_X \tOmega _{kl} \tOmega^{lj} D_X \tOmega_{ji} \geq 0.
    \]
    Combining with \eqref{eqn:C2-3rd-diff} gives
    \begin{equation}
    \label{neq:C2-eqn-2diff}
        \sum_{i=0}^{n-1} \frac{D^2_X \tOmega_{2i 2i+1}}{\tOmega_{2i 2i+1}}
        \geq D^2_X f + D^2_X \log \text{Pf}(\Omega_{ij}).
    \end{equation}
    Write $J$ in local coordinates as
    \[
    J = J^l_{\bar k} d\bar{z}^k \otimes \p_{z^l} + J^{\bar l}_k dz^k \otimes \p_{\bar{z}^l}.
    \]
    Notice that 
    \[
    \tOmega_{ij} = (\Omega_h)_{ij} + 
    \frac{1}{n-1}(S_1(\p \p_J u)\Omega_{ij} - 
    (-u_{i \overline k}J^{\overline k}_j + u_{j \overline k}J^{\overline k}_i ) ).
    \]
    Differentiating twice we get 
    \begin{equation}
    \begin{split}
       (n-1) D^2_X \tOmega_{2i2i+1} = &\,(n-1) D^2_X(\Omega_h)_{2i 2i+1} + 
     \sum_{p=0}^{2n-1} D^2_X u_{p \overline{p}} + S_1(\p \p_J u)D_X^2 \Omega_{2i2i+1} \\
    &- (D^2_X u_{2i \overline{2i}} + D^2_X u_{2i+1 \overline{2i+1}})
    + u_{i \overline{k}} D^2_X J^{\overline{k}}_j 
    - u_{j \overline{k}} D^2_X J^{\overline{k}}_i.
    \end{split}
    \end{equation}
    Here we used Remark \ref{rm:normal} again, namely at the point $p$,
    \[
    J^l_{\bar k,i} =J^{\bar l}_{k,i} =J^{\bar l}_{k,\bar i}=J^l_{\bar k, \bar i} =0.
    \]
    Combine with \eqref{neq:C2-eqn-2diff} 
    \begin{equation}
        \sum_{p=0}^{2n-1} D^2_X u_{p \overline{p}} \geq -C_5(\Tilde{\lambda} + 1).
    \end{equation}
    Then combining with \eqref{neq:C2-first-term-local} we obtain
    \begin{equation}
    \label{neq:C2-1-est}
        \frac{1}{n} \sum_{i=0}^{n-1} \sum_{p \neq i} 
        \frac{\Tilde{\lambda}_{2p\overline{2p}} 
        + \Tilde{\lambda}_{2p+1\overline{2p+1}}}{\tOmega_{2i2i+1}}
        \geq 
        -C'(\Tilde{\lambda} + 1).
    \end{equation}
    
    Now we have the eatimate of \eqref{eqn:C2-first-term}. The second term of \eqref{neq:C2-two-terms} has been dealt with in $C^1$ estimate as in \eqref{eqn:C1-1-local}
    \begin{equation}
    \begin{split}
       &\frac{\frac{1}{4}\p \p_J |du|_g^2\wedge A \wedge \bOmega^n}{ \tOmega{\mathstrut}^n \wedge \bOmega^n}\\
       = &- \frac{1}{n} \sum_{i = 0}^{n-1} \sum_{j = 0}^{2n-1} 
    \frac{(\Omega_h)_{2i2i+1, \overline{j}} u_{j} + (\bOmega_h)_{2i2i+1, j} u_{\overline{j}} }
    {\tOmega_{2i2i+1}} 
    + \frac{1}{n} \sum_{j=0}^{2n-1} 
    \frac{u_j(e^f)_{\overline{j}} + u_{\overline{j}}(e^f)_{j} }{e^f} \\
    &+ \frac{1}{n} \sum_{k=0}^{n-1} \sum_{j=0}^{2n-1} \sum_{i\neq k}
    \frac{|u_{2kj}|^2 + |u_{2k+1j}|^2}{\tOmega_{2i2i+1}} 
    + \frac{1}{n} \sum_{k=0}^{n-1} \sum_{j=0}^{2n-1} \sum_{i\neq k}
    \frac{|u_{2k\overline{j}}|^2 + |u_{2k+1\overline{j}}|^2}{\tOmega_{2i2i+1}} .
    \end{split} 
    \end{equation}
    Combining with \eqref{neq:C2-tOmega-bound} and $C^1$ estimate we obtain
    \begin{equation}
    \label{neq:C2-2-1}
    \frac{\frac{1}{4}\p \p_J |du|_g^2\wedge A \wedge \bOmega^n}{ \tOmega{\mathstrut}^n \wedge \bOmega^n}
    \geq -C_6 + C_7(|u_{ij}|^2 + |u_{i \bar j}|^2).
    \end{equation}
    By the definition of $\lambda$
    \begin{equation}
    \label{neq:C2-2-2}
    |u_{ij}|^2 + |u_{i \bar j}|^2 \geq C_8 \Tilde{\lambda}^2.
    \end{equation}
    Combining \eqref{neq:C2-2-1} and \eqref{neq:C2-2-2} we get
    \begin{equation}
    \label{neq:C2-2-est}
    \frac{\frac{1}{4}\p \p_J |du|_g^2\wedge A \wedge \bOmega^n}{ \tOmega{\mathstrut}^n \wedge \bOmega^n}
    \geq -C_6 + C\Tilde{\lambda}^2.
    \end{equation}
    Insert \eqref{neq:C2-1-est} and \eqref{neq:C2-2-est} into \eqref{neq:C2-two-terms}
    \begin{equation}
        0 \geq C\Tilde{\lambda}^2 - C' \Tilde{\lambda} - C''.
    \end{equation}
    This gives upper bound of $\Tilde{\lambda}$, therefore $\lambda$ is bounded above.
\end{proof}
\section{Proof of the Main Theorem}
\label{pmt}
Once we have the $C^2$ estimates,  the $C^{2,\alpha}$-estimates can be derived. In order to prove the main theorem, We consider the following continuity equation $(u_t, b_t)$ with $t\in[0,1]$:
\begin{align}
(\Omega_h+\frac{1}{n-1}(S_1(\p\p_J u_t)\Omega-\p\p_J u_t))^n=e^{tf+(1-t)f_0+b_t}\Omega^n,\label{main1}\\
\Omega_h+\frac{1}{n-1}(S_1(\p\p_J u_t)\Omega-\p\p_J u_t)>0,\quad \sup\limits_M u_t=0.\label{main2}
\end{align}
where $f_0=\log (\Omega_h^n/\Omega^n)$.
Consider the set
\[
S =\{
t \in [0, 1]:
  (u_t, b_t) \in C^{2,\alpha}(M, \mathbb{R}) \times \mathbb{R} 
\text{\ solves the equation\ }  \eqref{main1}, \eqref{main2} \}
\]

Clearly we have $0 \in S$. The $C^{2,\alpha}$-estimates implies closedness of $S$. We would like to show the openness as in \cite{blocki2013complex}.
Denote 
\[
\tOmega_{u} = \Omega_h+\frac{1}{n-1}(S_1(\p\p_J u)\Omega-\p\p_J u).
\]
Consider the operator 
\[
\mathcal{M} : \mathcal{A} \ni u \mapsto \frac{\tOmega_u^n}{\Omega_u} \in \mathcal{B},
\]
where
\begin{align*}
    &\mathcal{A} \coloneqq
    \{u \in C^{k+2, \alpha}(M): \tOmega_u > 0, \int_M u \Omega^n \wedge \bOmega^n = 0\}\\
    &\mathcal{B} \coloneqq
    \{\Tilde{f} \in C^{k, \alpha}(M): 
    \int_M \Tilde{f} \Omega^n \wedge \bOmega^n = \int_M \Omega^n \wedge \bOmega^n \}.
\end{align*}
It remains to show that for every $u \in \mathcal{A}$ 
the differential $d_{u}\mathcal{M}$ is an isomorphism. Indeed for $v \in T_u\mathcal{A}$ we have
\begin{align*}
    d_{u}\mathcal{M}(v) &= \frac{d}{dt}\bigg|_{t=0} \mathcal{M}(u + tv) 
    = \frac{d}{dt}\bigg|_{t=0} 
      \frac{(\Omega_h+\frac{1}{n-1}(S_1(\p\p_J(u+tv))\Omega-\p\p_J (u+tv)))^n}{\Omega^n}\\
    &= \frac{n}{n-1}
       \frac{(S_1(\p \p_J v)\Omega - \p \p_J v)\wedge \tOmega_u^{n-1}}{\Omega^n}
    = \frac{S_{n-1}(\tOmega_u) - 1}{2(n-1)}\Delta_{I, g}v.
\end{align*}
From general elliptic theory we know that the laplacian is a bijection between the space of functions of zero integral on $M$. Thus $\mathcal{M}$ is locally invertible and therefore $S$ is open.
\section*{}{\bf Acknowledgements:}  Fu is supported by  NSFC grant No. 12141104. Zhang is supported by NSFC grant  No. 11901102. 
\bibliographystyle{acm}
\bibliography{refs}

\end{document}